\documentclass{amsart}

\usepackage{amssymb, amsmath, amsthm}
\usepackage{mathrsfs}
\usepackage[fleqn,tbtags]{mathtools}
\usepackage[shortlabels]{enumitem}
\usepackage{tikz}
\usetikzlibrary{arrows.meta}
\usepackage{xcolor}
\usepackage[colorlinks,linkcolor={blue},citecolor={blue},urlcolor={blue!50!black}]{hyperref}

\let\le\leqslant
\let\leq\leqslant
\let\ge\geqslant
\let\geq\geqslant
\let\emptyset\varnothing

\theoremstyle{plain}
\newtheorem{theorem}{Theorem}[section]
\theoremstyle{remark}

\newtheorem{definition}{Definition}[section]
\theoremstyle{plain}

\newtheorem{corollary}[theorem]{Corollary}
\newtheorem{lemma}[theorem]{Lemma}
\numberwithin{equation}{section}
\newcommand{\D}{\{-1,1\}}
\newcommand{\E}{\mathbb E}

\title[Stein's inequality and superreflexivity]{Bourgain's vector-valued Stein inequality implies superreflexivity}

\author{Jan van Neerven}

\address{Delft Institute of Applied Mathematics\\
Delft University of Technology\\P.O. Box 5031\\2600 GA Delft\\The Netherlands}
\email{J.M.A.M.vanNeerven@tudelft.nl}

\begin{document}

\begin{abstract}
We prove that if $X$ is a Banach space for which Bourgain's vector-valued Stein inequality holds, then $X$ is superreflexive. We also present two examples of a superreflexive Banach space for which the vector-valued Stein property fails, the second of these being a twisted Hilbert space. This second example is used to deduce that the validity of the Stein inequality
is not a three-space property.
\end{abstract}



\maketitle
\enlargethispage{2pt}

\section{Introduction}

This paper continues the exploration, begun in \cite{Neerven-Stein}, of Bourgain's vector-valued Stein inequality in the context of vector-valued analysis.
The main result of the paper asserts that every Banach space $X$ for which Bourgain's vector-valued Stein inequality holds is superreflexive.
We present a relatively simple example of a superreflexive Banach space without the Stein property, and a more elaborate construction of a twisted Hilbert space with these properties.  The latter example implies that the validity of the Stein inequality
is not a three-space property.

Let us now explain the vector-valued Stein inequality in more detail.
Let $\mathscr F_k^{(d)}$ be the $\sigma$-algebra of $\D^d$ generated by the first
$k$ coordinate maps, with the understanding that $\mathscr F_0^{(d)} = \{\emptyset, \D^d\}$, and let $\mathcal E_k^{(d)}$ denote the corresponding conditional
expectation. We say that a Banach space $X$ \emph{satisfies Bourgain's
vector-valued Stein inequality}, or simply \emph{has the Stein property}, if there exists
a finite constant $C_X\ge 0$ such that for all dimensions $d\ge 1$ and integers $n\geq1$ we have
\begin{equation}\label{eq:PX-finite-Stein}
 \Bigl(\E_r\Bigl\|\sum_{j=1}^n r_j
                \mathcal E_{k_j}^{(d)}f_j\Bigr\|_{L^2(\D^d;X)}^2
       \Bigr)^{1/2}
 \leq C_X
 \Bigl(\E_r\Bigl\|\sum_{j=1}^n r_jf_j
                       \Bigr\|_{L^2(\D^d;X)}^2\Bigr)^{1/2}
\end{equation}
for all $k_1,\ldots,k_n\in\{0,\ldots,d\}$
and $f_1,\ldots,f_n\in L^2(\D^d;X)$. Here $\E_r$ denotes
expectation with respect to the \emph{Rademacher sequence} $(r_j)_{j\geq1}$
on a probability space $\Omega_r$. The variables $r_j$ are independent
and take the values $-1$ and $1$ with probability $1/2$ each.
The individual inequalities in \eqref{eq:PX-finite-Stein}, one for each
choice of $d$, $n$, and the $k_j$ and $f_j$, are referred to as the
{\em tests}. In this context the functions $f_j$ are
called the {\em inputs} and their conditional expectations
$\mathcal E_{k_j}^{(d)}f_j$ the {\em outputs}.

If $X$ has the Stein property, the least admissible constant in
\eqref{eq:PX-finite-Stein} is denoted by $\mathsf S(X)$
and is called the {\em (Rademacher) Stein constant} of $X$.

The classical scalar inequality for conditional expectations is
due to Stein \cite{Stein-Littlewood-Paley}. Bourgain established
its vector-valued extension for UMD spaces \cite{Bourgain-Stein}.
In particular, every UMD space has the Stein property.

By \cite[Theorems~3.1 and~3.2]{Neerven-Stein}, the Stein property is equivalent to the corresponding inequality on arbitrary probability spaces with
finite filtrations, with the inner $L^2$-norms in
\eqref{eq:PX-finite-Stein} replaced by $L^p$-norms for any fixed
$1<p<\infty$.  The Stein property is self-dual \cite{Eskenazis}, see also \cite[Proposition~3.5]{Neerven-Stein}: a Banach space $X$ has the Stein property if and
only if its dual $X^*$ does.

Over the complex scalars, one may replace the Rademacher
sequence by a \emph{Steinhaus sequence}, that is, a sequence of independent
random variables uniformly distributed on the unit circle. The best
constants for the two randomisations are comparable up to a universal
factor.

\section[Stein implies superreflexivity]{Stein implies superreflexivity}\label{sec:Stein-superreflexive}

A Banach space $Y$ is \emph{finitely
representable} in a Banach space $X$ if, for every
finite-dimensional subspace $F\subseteq Y$ and every
$\delta>0$, there is a linear isomorphism
$T:F\to T(F)\subseteq X$ whose {\em distortion} satisfies
$$\|T\|\,\|T^{-1}\|\leq1+\delta.$$
The space $X$ is \emph{superreflexive} if every Banach space
finitely representable in $X$ is reflexive.

\begin{theorem}
\label{thm:Stein-superreflexive}
Every Banach space with the Stein property is superreflexive.
\end{theorem}

Since superreflexivity implies $K$-convexity,
Theorem~\ref{thm:Stein-superreflexive} strengthens the result that
the Stein property implies $K$-convexity \cite[Lemma~4]{Eskenazis}.
It also removes the $\mathrm{UMD}^-$ assumption in
\cite[Lemma~6]{Eskenazis}.

It suffices to prove the theorem for real Banach spaces $X$. Indeed, the Stein
inequality for a complex space is the same inequality for
its underlying real space. If a complex Banach space
$Y$ is finitely representable in $X$, then its underlying real space
is finitely representable in the underlying real space of $X$.
Indeed, any finite-dimensional real subspace of $Y$ is contained
in its finite-dimensional complex span, and a complex embedding
of the span restricts to a real embedding of the subspace with
at most the same distortion. Since real and complex reflexivity
are equivalent for $Y$, superreflexivity of the underlying real
space of $X$ implies superreflexivity of $X$. We therefore work
over the real scalars in the rest of this section.

\subsection{Reduction to a basis with equal-sign additivity}

We first record two elementary stability properties of the Stein inequality.
In their proofs we write $L^2(Z)$ for
$L^2(\Omega_r\times\D^d;Z)$ with the product probability measure,
where $\D^d$ is the discrete cube of the test under consideration.

\begin{lemma}\label{lem:sr-finite-representability}
If $X$ has the Stein property and $Y$ is finitely representable in $X$, then
 $Y$ has the Stein property and
$$\mathsf S(Y)\leq\mathsf S(X).$$
\end{lemma}
\begin{proof} Consider a test in
\eqref{eq:PX-finite-Stein} with values in $Y$. The values of its finitely
many input functions span a finite-dimensional space $F$.
For each $\delta>0$, finite representability gives an
embedding $T:F\to X$ with distortion at most $1+\delta$. The
conditional averages of the $f_j$ belong to $F$, and taking these
averages commutes with applying $T$. The Stein inequality for the
functions $Tf_j$ and the bound for $T^{-1}$ on $T(F)$ give
\begin{align*}
 \Bigl\|\sum_jr_j\mathcal E_{k_j}^{(d)}f_j\Bigr\|_{L^2(Y)}
 &\leq\|T^{-1}\|\mathsf S(X)
       \Bigl\|\sum_jr_jTf_j\Bigr\|_{L^2(X)}\\
 &\leq(1+\delta)\mathsf S(X)
       \Bigl\|\sum_jr_jf_j\Bigr\|_{L^2(Y)}.
\end{align*}
Letting $\delta\downarrow0$ gives the required bound.
\end{proof}

\begin{lemma}\label{lem:sr-ell2-sum}
If $X$ has the Stein property, then the direct sum $X\oplus_2 X$, with norm
$\|(x,y)\|=(\|x\|^2+\|y\|^2)^{1/2}$, has the Stein property, and
$$\mathsf S(X\oplus_2 X) \le \mathsf S(X).$$
\end{lemma}
\begin{proof} If $f_j=(g_j,h_j)$, then
\[
 \Bigl\|\sum_jr_jf_j\Bigr\|_{L^2(X\oplus_2X)}^2
 =\Bigl\|\sum_jr_jg_j\Bigr\|_{L^2(X)}^2
  +\Bigl\|\sum_jr_jh_j\Bigr\|_{L^2(X)}^2.
\]
The same identity holds for the conditional averages. Adding
the two squared Stein estimates therefore gives the assertion.
\end{proof}

We write $c_{00}$ for the real sequences with finite support.
A \emph{Schauder basis} $(e_i)$ of a Banach space is a sequence
such that every vector has a unique norm-convergent expansion
$\sum_i a_i e_i$. It is \emph{normalised} if $\|e_i\|=1$ for
every $i$.

\begin{definition} A Schauder basis is \emph{equal-signs-additive}, abbreviated
ESA, if it is normalised and satisfies
\begin{equation}
\label{eq:sr-ESA}
 \Bigl\|\sum_{i<k}a_i e_i+(a_k+a_{k+1})e_k
                  +\sum_{i>k+1}a_i e_i\Bigr\|
 =\Bigl\|\sum_i a_i e_i\Bigr\|,\qquad
 a_k a_{k+1}\geq0,
\end{equation}
for every $a\in c_{00}$ and every $k$.
Thus adjacent coefficients of the same sign may be combined
without changing the norm.
\end{definition}

We use two standard consequences of this definition.
First, the norm is \emph{invariant under spreading}: for every
strictly increasing sequence of indices $(n_i)$,
\[
 \Bigl\|\sum_i a_i e_{n_i}\Bigr\|
 =\Bigl\|\sum_i a_i e_i\Bigr\|.
\]
Equivalently, inserting or deleting zero coefficients preserves
the norm. Secondly, the equality in \eqref{eq:sr-ESA} becomes
an inequality with left-hand side no larger than the right-hand
side when the sign condition $a_k a_{k+1}\geq0$ is omitted.
Both assertions follow from
\cite[Lemma~1]{BrunelSucheston1976}; see also
\cite[Lemma~2.2]{OstrovskiiRandrianantoanina-ESA}.
Applying the sign-free adjacent compression repeatedly within
each of a finite collection of consecutive blocks replaces
each block by its sum, with zeros in the vacated positions.
Deleting the vacated positions by spreading invariance gives the
contractive block-compression estimate
\begin{equation}
\label{eq:sr-compression}
 \Bigl\|\sum_{j=1}^d
       \Bigl(\sum_{i=k_{j-1}+1}^{k_j}a_i\Bigr)e_j\Bigr\|
 \leq\Bigl\|\sum_{i=1}^N a_i e_i\Bigr\|,
 \qquad 0=k_0<\cdots<k_d=N.
\end{equation}

We use the following consequence of Brunel and Sucheston's work
\cite{BrunelSucheston1975}, in the formulation of
\cite[Theorem~2.3]{OstrovskiiRandrianantoanina-ESA}. Taking $a_{k+1}=1$
and all other coefficients zero in \eqref{eq:sr-ESA} gives
$\|e_k\|=\|e_{k+1}\|$, so a single rescaling normalises the basis
in that result and preserves equal-sign additivity.

\begin{lemma}\label{lem:sr-brunel-sucheston}
If $Y$ is a nonreflexive real Banach space, then there
exists a real Banach space $E_0$ which is finitely representable in $Y$
and has an ESA Schauder basis.
\end{lemma}

\begin{lemma}
\label{lem:sr-reduction}
If $X$ has the Stein property and is not superreflexive,
there exists a real
Stein space $E$ with a normalised ESA basis $(e_i)$ such that
\[
 \Bigl\|\sum_{i=1}^N a_i e_i\Bigr\|_E
 =\Bigl\|\sum_{i=1}^N a_{N+1-i}e_i\Bigr\|_E,
 \qquad a_1,\ldots,a_N\in\mathbb R.
\]
\end{lemma}

\begin{proof}
By the reduction to real scalars above, we may assume that $X$ is real.
Since $X$ is not superreflexive, there exists a nonreflexive real
Banach space $Y$ finitely representable in $X$. Lemma~\ref{lem:sr-finite-representability}
gives the Stein property for $Y$, with $\mathsf S(Y)\leq\mathsf S(X)$.
Lemma~\ref{lem:sr-brunel-sucheston}, applied to $Y$, gives a real Banach
space $E_0$, finitely representable in $Y$, with an ESA Schauder basis
$(e_i^0)_{i\geq1}$. A second application of
Lemma~\ref{lem:sr-finite-representability} gives the Stein property
for $E_0$, with $\mathsf S(E_0)\leq\mathsf S(Y)$.
By Lemma~\ref{lem:sr-ell2-sum}, the Hilbertian sum
$Z:=E_0\oplus_2 E_0$ has the Stein property and
$\mathsf S(Z)\leq\mathsf S(E_0)$.

The interval construction below uses the coefficients in both their
original order and the reverse order. We therefore give $c_{00}$ a norm
that treats the two orders symmetrically. For $N\geq1$, let $c_{00}^{(N)}$ consist
of the sequences in $c_{00}$ supported in $\{1,\ldots,N\}$.
For $a\in c_{00}^{(N)}$ we define $\Phi_N(a)\in E_0$ and $R_Na \in c_{00}^{(N)}$ by
\[
 \Phi_N(a):=\sum_{i=1}^N a_i e_i^0,
 \qquad
 R_Na:=(a_N,\ldots,a_1,0,0,\ldots).
\]
Thus $R_N$ is a linear operator on $c_{00}^{(N)}$ and $R_N^2$ is the identity mapping on $c_{00}^{(N)}$.
The mapping $J_N:c_{00}^{(N)}\to Z$ given by
\[
 J_Na:=\frac1{\sqrt{2}}\bigl(\Phi_N(a),\Phi_N(R_Na)\bigr)
\]
is linear, and it is injective because $\Phi_N(a)=0$ implies
$a_1=\cdots=a_N=0$ by linear independence of the basis vectors.
In particular, the formula
\begin{equation}\label{eq:sr-reversal}
 |a|_E^2
 :=\frac12\|\Phi_N(a)\|_{E_0}^2
   +\frac12\|\Phi_N(R_Na)\|_{E_0}^2
 =\|J_Na\|_Z^2,
 \qquad a\in c_{00}^{(N)},
\end{equation}
defines a norm on $c_{00}^{(N)}$.

For $a\in c_{00}^{(N)}$ and $M\geq N$, the entries of $a$ after the $N$-th
are zero, so $\Phi_M(a)=\Phi_N(a)$. The first $M$ entries of $R_Ma \in c_{00}^{(M)}$
consist of $M-N$ zeros followed by $a_N,\ldots,a_1$. Hence
\[
 \Phi_M(R_Ma)=\sum_{i=1}^N a_{N+1-i}e_{M-N+i}^0.
\]
By spreading invariance of $(e_i^0)$, we have
\[
 \|\Phi_M(R_Ma)\|_{E_0}
 =\Bigl\|\sum_{i=1}^N a_{N+1-i}e_i^0\Bigr\|_{E_0}
 =\|\Phi_N(R_Na)\|_{E_0}.
\]
Thus \eqref{eq:sr-reversal} assigns the same norm to $a$ whether it
is evaluated in $c_{00}^{(N)}$ or in $c_{00}^{(M)}$. Every pair of finitely supported
sequences and their sum lie in a common $c_{00}^{(N)}$, so the triangle
inequality on the $c_{00}^{(N)}$ gives the triangle inequality on $c_{00}$.
Definiteness and homogeneity also follow from the norms on the $c_{00}^{(N)}$.
We let $E$ be the completion of $c_{00}$ in this norm and denote the
canonical unit vectors in $E$ by $(e_i)_{i\geq1}$. We identify $c_{00}^{(N)}$
with $\operatorname{span}\{e_1,\ldots,e_N\}$ in $E$. By
\eqref{eq:sr-reversal}, $J_N$ is an isometric embedding of this
subspace into $Z$.

For $a\in c_{00}^{(N)}$, replacing $a$ by $R_Na$ in \eqref{eq:sr-reversal}
interchanges the two terms. Thus
\[|R_Na|_E=|a|_E.\]
Taking $a=e_i$ and $N=i$ in \eqref{eq:sr-reversal} gives
\[
 \|e_i\|_E^2=\frac12\|e_i^0\|_{E_0}^2+\frac12\|e_1^0\|_{E_0}^2=1,
\]
where we used the normalisation of the basis of $E_0$.

We next verify that $(e_i)_{i\geq1}$ is a Schauder basis in $E$. Let $Q_m$ be the
$m$-th partial-sum projection of $(e_i^0)_{i\geq1}$ on $E_0$, with $Q_0=0$.
Since $(e_i^0)_{i\geq1}$ is a Schauder basis for $E_0$, the projections $Q_m$ are uniformly bounded on $E_0$, say with bound
\[K_0:=\sup_{m\geq0}\|Q_m\|<\infty.\]
Moreover, $Q_mx\to x$ for every $x\in E_0$.

For $m\geq1$, put $V_m:=\operatorname{span}\{e_1,\ldots,e_m\}$.
On $c_{00}$, let $P_m$ retain the first $m$ entries of a sequence
and set the later entries to zero. For $a\in c_{00}^{(N)}$ with $N\geq m$,
we have
\begin{align*}
 \Phi_N(P_ma)&=Q_m\Phi_N(a),\\
 \Phi_N(R_NP_ma)&=(I-Q_{N-m})\Phi_N(R_Na).
\end{align*}
In the second identity, the entries $a_{m+1},\ldots,a_N$ that
$P_m$ sets to zero occupy the first $N-m$ positions after reversal.
Since $\|I-Q_{N-m}\|\leq 1+K_0$, the two identities give
\begin{align*}
 |P_ma|_E^2
 \leq \frac12 K_0^2\|\Phi_N(a)\|_{E_0}^2
      +\frac12(1+K_0)^2\|\Phi_N(R_Na)\|_{E_0}^2
 \leq (1+K_0)^2|a|_E^2.
\end{align*}
Every $a\in c_{00}$ lies in some $c_{00}^{(N)}$ with $N\geq m$, so $P_m$
extends by continuity to a bounded operator on $E$, with
$\|P_m\|\leq1+K_0$. The extension takes values in $V_m$ because
$P_m(c_{00})\subseteq V_m$ and $V_m$ is closed in $E$. Since $P_m$
fixes every element of $V_m$, its range is $V_m$.

For $x\in E$ and $b\in c_{00}$, we have $P_mb=b$ for every
sufficiently large $m$. For these $m$,
\[
 \|P_mx-x\|_E
 \leq \|P_m(x-b)\|_E+\|x-b\|_E
 \leq (2+K_0)\|x-b\|_E.
\]
As $b$ can be chosen arbitrarily close to $x$, this proves
$P_mx\to x$ in norm. With $P_0=0$, the bounded operator
$P_m-P_{m-1}$ takes values in $\mathbb R e_m$, since it does so
on $c_{00}$ and the one-dimensional subspace $\mathbb R e_m$ is closed. For each $x\in E$,
let $\alpha_m(x)$ be the scalar such that
\[(P_m-P_{m-1})x=\alpha_m(x)e_m.\]
Summing these identities gives
$P_mx=\sum_{i=1}^m\alpha_i(x)e_i$, so the series
$\sum_{i\geq1}\alpha_i(x)e_i$ converges to $x$. If $\sum_{i\geq1}b_i e_i$
also converges to $x$, applying $P_m-P_{m-1}$ to its partial sums
of order at least $m$ and passing to the limit gives
$b_me_m=(P_m-P_{m-1})x$. Thus $b_m=\alpha_m(x)$ for every $m$,
which proves uniqueness of the expansion.

For equal-sign additivity, let $a\in c_{00}^{(N)}$ and $1\leq k<N$ satisfy
$a_ka_{k+1}\geq0$. Let $b\in c_{00}^{(N)}$ be
obtained by replacing $(a_k,a_{k+1})$ with $(a_k+a_{k+1},0)$.
The ESA identity in $E_0$ gives
\[\|\Phi_N(b)\|_{E_0}=\|\Phi_N(a)\|_{E_0}.\]
In the reversed sequences the two entries
occur at positions $N-k$ and $N-k+1$. At these positions $R_Na$
has the pair $(a_{k+1},a_k)$ and $R_Nb$ has the pair $(0,a_k+a_{k+1})$.
In $E_0$ the ESA identity allows each of these pairs to be replaced
by $(a_k+a_{k+1},0)$ without changing the norm, since $a_{k+1}a_k\geq0$
and $0\cdot(a_k+a_{k+1})=0$. All the other entries agree, so
\[\|\Phi_N(R_Nb)\|_{E_0}=\|\Phi_N(R_Na)\|_{E_0}.\] Consequently
$|b|_E=|a|_E$ by \eqref{eq:sr-reversal}. Any pair of adjacent
coefficients in a finitely supported sequence can be treated by
taking $N$ large enough. This proves \eqref{eq:sr-ESA} for the
normalised Schauder basis $(e_i)_{i\geq1}$ of $E$. In particular, the
spreading invariance recalled above applies to $(e_i)_{i\geq1}$.

We will now show that $E$ is finitely representable in $Z$. Let
$F$ be a nonzero finite-dimensional subspace of $E$ and fix
$0<\varepsilon<1$. By compactness of the unit sphere of $F$,
it contains a finite subset $S$ such that each unit vector
$x\in F$ lies within $\varepsilon/(2(2+K_0))$ of some point $y\in S$.
Since $P_Ny\to y$ for every $y\in S$, for all sufficiently large $N$ we have
$\|P_Ny-y\|_E<\varepsilon/2$ for every $y\in S$.
For a unit vector $x\in F$ and a corresponding $y\in S$ we then have
\[
 \|P_Nx-x\|_E \leq (2+K_0)\|x-y\|_E+\|P_Ny-y\|_E <\varepsilon = \varepsilon \|x\|.
\]
By homogeneity and the triangle inequality, it follows that
\[
 (1-\varepsilon)\|x\|_E\leq\|P_Nx\|_E
 \leq(1+\varepsilon)\|x\|_E, \qquad x\in F.
\]
The linear map $T:F\to Z$ given by $Tx=J_N(P_Nx)$ satisfies these
same bounds, since $J_N$ is an isometry on $c_{00}^{(N)}$. The lower bound
makes $T$ injective and gives $\|T^{-1}\|\leq (1-\varepsilon)^{-1}$
on its image. The upper bound gives $\|T\|\leq 1+\varepsilon$.
For any $\delta>0$, choosing $\varepsilon=\delta/(2+\delta)$ gives
\[
 \|T\|\,\|T^{-1}\|
 \leq \frac{1+\varepsilon}{1-\varepsilon}=1+\delta.
\]
This proves finite representability. Applying
Lemma~\ref{lem:sr-finite-representability} to $E$ and $Z$ gives
the Stein property for $E$, with
\[
 \mathsf S(E)\leq\mathsf S(Z)\leq\mathsf S(E_0)
 \leq\mathsf S(Y)\leq\mathsf S(X).
\]
\end{proof}

We now fix a space $E$ as in
Lemma~\ref{lem:sr-reduction}. We write $e_i^*$ for its coordinate
functionals, characterised by $e_i^*(\sum_k a_k e_k)=a_i$, and put
\begin{align}\label{eq:defP-K}
 P_m\Bigl(\sum_i a_i e_i\Bigr):=\sum_{i=1}^m a_i e_i,
 \qquad K:=\sup_m\|P_m\|<\infty.
\end{align}
For $a\in c_{00}$, let
\begin{align}\label{eq:def-rhoa}
 |a|_E:=\Bigl\|\sum_i a_i e_i\Bigr\|_E,
 \qquad
 \rho(a):=\Bigl(\E_r
              \Bigl\|\sum_i r_i a_i e_i\Bigr\|_E^2\Bigr)^{1/2}.
\end{align}
The norm $\rho$ is unchanged by arbitrary sign
changes, by insertion of zeros and by finite reversal.
For a set $J$ of coordinates, independence and the mean value property
$\E_r r_i=0$ of the Rademacher variables give
\[
 \E_r\Bigl(\sum_i r_i a_i e_i\,\Big|\,r_j:j\in J\Bigr)
 =\sum_{i\in J}r_i a_i e_i.
\]
Contractivity of this conditional expectation on $L^2(\Omega_r;E)$
shows that setting the other coordinates to zero cannot
increase $\rho$.

Since $a\in c_{00}$ has finite support, there is an $N\geq1$
such that $a_i=0$ for $i>N$. Applying
\eqref{eq:sr-compression} with a single block, that is, with
$d=1$, $k_0=0$ and $k_1=N$, gives
\[
 \Bigl\|\Bigl(\sum_{i=1}^N a_i\Bigr)e_1\Bigr\|_E
 \leq \Bigl\|\sum_{i=1}^N a_i e_i\Bigr\|_E
 =|a|_E.
\]
The basis is normalised, so the expression on the left equals
$\bigl|\sum_{i=1}^N a_i\bigr|\|e_1\|_E
=\bigl|\sum_i a_i\bigr|$. Hence
$\bigl|\sum_i a_i\bigr|\leq |a|_E$.
 Applying this scalar estimate to the coefficients $(r_ia_i)$
and using $\E_r(r_ir_j)=\mathbf1_{\{i=j\}}$ gives
\[
 \sum_i a_i^2=\E_r\Bigl|\sum_i r_i a_i\Bigr|^2
 \leq\E_r\Bigl\|\sum_i r_i a_i e_i\Bigr\|_E^2.
\]
Taking square roots gives
\begin{equation}
\label{eq:sr-lower}
 \|a\|_{\ell^2}:=\Bigl(\sum_i a_i^2\Bigr)^{1/2}
 \leq\rho(a).
\end{equation}

\subsection{Two applications of Stein's inequality}

Our next aim is to prove the estimates
\[
 \rho(a)\leq A|a|_E,
 \qquad
 \rho(a)\leq B\|a\|_{\ell^2},
 \qquad a\in c_{00},
\]
with constants independent of the number of coefficients.
The first estimate will follow from Stein's inequality
in $E$. For the second, we will apply Stein's inequality
in $E^*$ and then use finite-dimensional duality.

Both arguments use a random construction of intervals
of basis indices. The chosen intervals will be disjoint,
while their conditional averages will correspond to
nested parent intervals. Disjointness will allow us to
estimate the inputs in Stein's inequality. Nesting will
allow us to separate all left endpoints from all right
endpoints in the first argument, and to find a coordinate
common to all parent intervals in the second.

For a fixed
$n\geq1$, the construction takes place on the probability space $\Omega_n = \D^n$ with its
coordinate filtration $(\mathscr F_j)_{j=0}^n$.
We start with the integer interval $I_0=[1,1+8^n)$.
At step $j$, the $j$-th coordinate chooses, with equal probability,
one of the two halves $H_j$ of $I_{j-1}$. The next interval
$I_j$ is the middle half of the unchosen half.

More explicitly, if $I_{j-1}=[\lambda,\lambda+L)$, the two
possibilities are
\begin{align*}
\bullet \  H_j &= [\lambda,\lambda+L/2),   &&\hskip-1.5truecm I_j =  [\lambda+5L/8,\lambda+7L/8);\\
\bullet \  H_j &= [\lambda+L/2,\lambda+L), &&\hskip-1.5truecm I_j =  [\lambda+L/8,\lambda+3L/8).
\end{align*}

In either case $|I_j|=|I_{j-1}|/4$, so
$|I_{j-1}|=2^{3n-2j+2}$. For $1\leq j\leq n$ this
length is divisible by $8$. Thus every endpoint is an
integer, and the displayed gaps have positive integer length.
The choice at step $j$ depends only on the first $j$
coordinates, so $I_j$ and $H_j$ are $\mathscr F_j$-measurable.

The intervals $H_1,\ldots,H_n$ are disjoint and separated by
nonempty gaps. Their spatial order consists of the left-half
choices in increasing chronological order, followed by the
right-half choices in decreasing chronological order. This is
because every subsequent interval lies strictly to the right
of a chosen left half, and strictly to the left of a chosen
right half. The intervals $I_j$ are strictly nested at both
endpoints.
The reversal symmetry introduced in
Lemma~\ref{lem:sr-reduction} will allow us to handle
the right-half choices in this reversed order.

For an interval $H=[\lambda,\mu)$ with integer endpoints
satisfying $1\leq\lambda<\mu$, define
\begin{align}\label{eq:defvHhH}
 v_H:=e_\lambda-e_\mu\in E,
 \qquad
 h_H:=\sum_{k=\lambda}^{\mu-1}e_k^*\in E^*.
\end{align}
For an integer $m$ with $\lambda<m<\mu$ we have $v_{[\lambda,m)}+v_{[m,\mu)}=v_{[\lambda,\mu)}$ and
$h_{[\lambda,m)}+h_{[m,\mu)}=h_{[\lambda,\mu)}$.
Conditional on $\mathscr F_{j-1}$, the parent $I_{j-1}$
is fixed, and the independent $j$-th sign chooses each half
with probability $1/2$. Averaging these two values therefore
gives
\begin{equation}
\label{eq:sr-averages}
 \mathbb E(v_{H_j}\mid\mathscr F_{j-1})=\frac12v_{I_{j-1}},
 \qquad
 \mathbb E(h_{H_j}\mid\mathscr F_{j-1})=\frac12h_{I_{j-1}}.
\end{equation}

\begin{lemma}
\label{lem:sr-randomisation}
There is a constant $A<\infty$ such that
$\rho(a)\leq A|a|_E$ for every $a\in c_{00}$.
\end{lemma}

\begin{proof}
For $a=(a_1,\ldots,a_n)$ and $J\subseteq\{1,\ldots,n\}$, put
\[U(a):=\sum_{j=1}^n a_j(e_{2j-1}-e_{2j}), \qquad
U_J(a):=\sum_{j\in J}a_j(e_{2j-1}-e_{2j}).\]
The block-compression estimate \eqref{eq:sr-compression},
followed by spreading invariance, allows each omitted pair
$(a_j,-a_j)$ to be replaced by $(0,0)$. Hence
$\|U_J(a)\|_E\leq\|U(a)\|_E$. If $\sigma_j\in\{-1,1\}$
and $J_\pm=\{j\in J:\sigma_j=\pm1\}$, then
\[
 \Bigl\|\sum_{j\in J}\sigma_j a_j(e_{2j-1}-e_{2j})\Bigr\|_E
 =\|U_{J_+}(a)-U_{J_-}(a)\|_E
 \leq2\|U(a)\|_E.
\]
Reversing all $2n$ coordinates of $U(a)$ gives
$-U(a_n,\ldots,a_1)$, so reversal invariance also permits
the pairs to be read in the opposite order. Finally,
spreading invariance applied to the odd and even coordinates
separately gives
\[
 \|U(a)\|_E
 \leq\Bigl\|\sum_j a_j e_{2j-1}\Bigr\|_E
     +\Bigl\|\sum_j a_j e_{2j}\Bigr\|_E
 =2|a|_E.
\]

We apply \eqref{eq:PX-finite-Stein} to $f_j=a_jv_{H_j}$ with
$k_j=j-1$. For a fixed $\omega\in\Omega_n$, let $L$ and $R$ be the sets of
indices at which the left and right halves are chosen.
The endpoints of the intervals indexed by $L$ occur in
chronological order, whereas those indexed by $R$ occur
in reverse chronological order. Within either group,
spreading invariance identifies the sum of endpoint vectors
with the corresponding signed subsequence of pairs in
$U(a)$, after reversal for the second group.
Each group therefore has norm at most $2\|U(a)\|_E$.
The triangle inequality gives, for every choice of the
external signs, the pointwise bound
\begin{align}\label{eq:sr-input-bound}
 \Bigl\|\sum_{j=1}^n r_j a_j v_{H_j}\Bigr\|_E
 \leq4\|U(a)\|_E\leq8|a|_E.
\end{align}

Writing $I_{j-1}=[\lambda_{j-1},\mu_{j-1})$, strict nesting gives
$$\lambda_0<\cdots<\lambda_{n-1}<\mu_{n-1}<\cdots<\mu_0.$$
By \eqref{eq:sr-averages}, the randomised output at
$\omega\in\Omega_n$ is
$\frac12\sum_jr_ja_jv_{I_{j-1}(\omega)}$.
Every left endpoint is at most $\lambda_{n-1}(\omega)$, while
every right endpoint is strictly larger. Hence the definition
\eqref{eq:defP-K} of the partial-sum projections gives
\[
 P_{\lambda_{n-1}(\omega)}\Bigl(\frac12\sum_{j=1}^n r_j a_j
                         v_{I_{j-1}(\omega)}\Bigr)
 =\frac12\sum_{j=1}^n r_j a_j e_{\lambda_{j-1}(\omega)}.
\]

The indices
$\lambda_0(\omega),\ldots,\lambda_{n-1}(\omega)$ are strictly
increasing. Spreading invariance and the norm estimate $\|P_m\|\leq K$ therefore give

\begin{equation}\label{eq:sr-output-bound}
\begin{aligned}
 \frac{1}{4}\rho(a)^2
 & = \E_r
              \Bigl\|\frac12\sum_i r_i a_i e_i\Bigr\|_E^2
 \\ & =\E_r\Bigl\|\frac12\sum_jr_ja_je_{\lambda_{j-1}(\omega)}
       \Bigr\|_E^2
 \leq K^2\E_r\Bigl\|\frac12\sum_jr_ja_j
       v_{I_{j-1}(\omega)}\Bigr\|_E^2.
\end{aligned}
\end{equation}

This estimate
holds at every \(\omega\in\Omega_n\) with the same constant $K$.
Integrating \eqref{eq:sr-output-bound} over $\Omega_n$ gives
\[
 \frac14\rho(a)^2
 \leq K^2\E_\omega\E_r
 \Bigl\|\frac12\sum_{j=1}^n r_j a_j
                    v_{I_{j-1}(\omega)}\Bigr\|_E^2.
\]
By \eqref{eq:sr-averages}, the randomised sum on the right
is the output corresponding to the inputs
$f_j=a_jv_{H_j}$ and the conditional expectations
$\mathcal E_{j-1}^{(n)}$. Taking square roots therefore gives
\[
 \frac{\rho(a)}{2K}
 \leq
 \Bigl\|\sum_{j=1}^n r_j\mathcal E_{j-1}^{(n)}f_j
       \Bigr\|_{L^2(\Omega_n\times\Omega_r;E)}.
\]
On the other hand, \eqref{eq:sr-input-bound} holds for every
$\omega\in\Omega_n$ and every choice of the external signs.
Squaring and integrating this bound gives
\[
 \Bigl\|\sum_{j=1}^n r_jf_j
       \Bigr\|_{L^2(\Omega_n\times\Omega_r;E)}
 \leq8|a|_E.
\]
Applying \eqref{eq:PX-finite-Stein} with the Stein constant
$\mathsf S(E)$, we conclude that
\begin{align*}
 \frac{\rho(a)}{2K}
 \leq
 \Bigl\|\sum_{j=1}^n r_j\mathcal E_{j-1}^{(n)}f_j
       \Bigr\|_{L^2(\Omega_n\times\Omega_r;E)}
 \leq \mathsf S(E)
 \Bigl\|\sum_{j=1}^n r_jf_j
       \Bigr\|_{L^2(\Omega_n\times\Omega_r;E)}
 \leq8\mathsf S(E)|a|_E.
\end{align*}
Thus
\[
 \rho(a)\leq16K\mathsf S(E)|a|_E,
\]
which proves the assertion with $A:=16K\mathsf S(E)$.
\end{proof}

For $b\in c_{00}$, define the dual norm
\begin{equation}\label{eq:sr-dual-norm}
 \rho^*(b):=\sup\Bigl\{
       \Bigl|\sum_i a_i b_i\Bigr|:\,
       a\in c_{00},\ \rho(a)\leq1\Bigr\}.
\end{equation}
By \eqref{eq:sr-lower} and the Cauchy--Schwarz inequality,
\[
 \Bigl|\sum_i a_i b_i\Bigr|
 \leq\|a\|_{\ell^2}\|b\|_{\ell^2}
 \leq\rho(a)\|b\|_{\ell^2},
\]
so
\begin{align}\label{eq:rhostarb}
\rho^*(b)\leq\|b\|_{\ell^2}<\infty.
\end{align}

If $b\in c_{00}$ is nonzero, say $b_i\ne0$, testing with the coefficient
vector having a single entry $1$ at position $i$ gives
$\rho^*(b)\geq|b_i|\ne 0$, since $\|e_i\|_E=1$.
Scalar homogeneity and the triangle inequality being clear,
we see that $\rho^*$ is a norm on $c_{00}$.

For a finite set $J$ containing the support of $b$, let $a_J$
denote the restriction of $a$ to $J$, with all other entries zero.
We claim that the supremum defining $\rho^*(b)$ can be restricted
to such vectors. More precisely,
\begin{equation}\label{eq:sr-dual-restriction}
 \rho^*(b)=
 \sup\Bigl\{\Bigl|\sum_{i\in J}a_i b_i\Bigr|:\,
       \operatorname{supp}a\subseteq J,\ \rho(a)\leq1\Bigr\}.
\end{equation}
The right-hand side is at most $\rho^*(b)$, since every vector
admissible there is also admissible in \eqref{eq:sr-dual-norm}.
For the opposite inequality, take any $a\in c_{00}$ with
$\rho(a)\leq1$. The suppression property proved after
\eqref{eq:def-rhoa} gives
\[
 \rho(a_J)\leq\rho(a)\leq1.
\]
Also, since $b_i=0$ outside $J$,
\[
 \sum_i a_i b_i=\sum_{i\in J}(a_J)_i b_i.
\]
Thus $a_J$ is admissible in the restricted supremum and gives
exactly the same value of the pairing as $a$.
Every value occurring in the full supremum is therefore bounded
by the restricted supremum. Taking the supremum over $a$ proves
the opposite inequality and hence \eqref{eq:sr-dual-restriction}.

We identify $\mathbb R^J$ with the sequences supported in $J$ by
extending their entries by zero. The right-hand side of
\eqref{eq:sr-dual-restriction} is precisely the operator norm of
the functional $c\mapsto\sum_{i\in J}c_i b_i$ on this
finite-dimensional space equipped with $\rho$.
Thus its finite-dimensional dual norm equals the globally
defined norm $\rho^*(b)$.

The norm $\rho^*$ enjoys the same sign invariance, reversal,
and spreading symmetries as $\rho$.
For sign invariance, a change of variables $a_i\mapsto\sigma_i a_i$,
where $\sigma_i\in\{-1,1\}$, preserves the constraint $\rho(a)\leq1$
and moves the signs from $b$ to $a$ in the pairing.
For reversal, \eqref{eq:sr-dual-restriction} allows us
to restrict to a finite coordinate space containing the support of $b$.
Reversing both coefficient lists
preserves their pairing, and reversal invariance of $\rho$ preserves
the set over which the supremum is taken.
For spreading, suppose that the entries $b_1,\ldots,b_d$ are placed
at increasing positions $n_1<\cdots<n_d$. In computing the dual norm
of this new vector, \eqref{eq:sr-dual-restriction} allows us
to retain only the coordinates $n_1,\ldots,n_d$ of the testing vector.
If $a'$ is this restriction and $c_i:=a_{n_i}$, spreading invariance
gives $\rho(a')=\rho((c_i)_{i=1}^d)$, while the pairing becomes
$\sum_{i=1}^d c_i b_i$. This substitution identifies the two
suprema, proving spreading invariance for $\rho^*$.
Finally, for any set $J$ and the corresponding restriction $b_J$,
\begin{equation}\label{eq:sr-dual-suppression}
 \rho^*(b_J)
 =\sup_{\rho(a)\leq1}\Bigl|\sum_i(a_J)_i b_i\Bigr|
 \leq\rho^*(b),
\end{equation}
because $\rho(a_J)\leq\rho(a)\leq1$.

Finite-dimensional duality also gives, for every finite $J$,
\begin{equation}\label{eq:sr-coefficient-duality}
 \rho(a)=
 \sup\Bigl\{\Bigl|\sum_i a_i b_i\Bigr|:\,
       \operatorname{supp}b\subseteq J,\ \rho^*(b)\leq1\Bigr\},
 \qquad\operatorname{supp}a\subseteq J.
\end{equation}
One inequality follows from the definition
of $\rho^*$. For the reverse inequality when $a\ne0$, the Hahn--Banach
theorem supplies a functional of norm one on $(\mathbb R^J,\rho)$
whose value at $a$ is $\rho(a)$. Writing that functional as
$c\mapsto\sum_{i\in J}c_i b_i$ gives the required testing vector $b$.
Equation \eqref{eq:sr-dual-restriction} identifies its finite-dimensional
dual norm with $\rho^*(b)$. The case $a=0$ is immediate.

\begin{lemma}
\label{lem:sr-euclidean}
There is a constant $B<\infty$ such that
$\rho(a)\leq B\|a\|_{\ell^2}$ for every $a\in c_{00}$.
\end{lemma}
\begin{proof}
We will first prove the dual estimate
\begin{equation}\label{eq:sr-dual-euclidean}
 \|b\|_{\ell^2}\leq B\rho^*(b),
 \qquad b\in c_{00},
\end{equation}
which reverses \eqref{eq:rhostarb} up to a constant.
Once this is known, \eqref{eq:sr-coefficient-duality}, the
Cauchy--Schwarz inequality, and \eqref{eq:sr-dual-euclidean} give
\begin{align}\label{eq:rho-Brhostar}
 \rho(a)
 \leq\|a\|_{\ell^2}\sup_{\rho^*(b)\leq1}\|b\|_{\ell^2}
 \leq B\|a\|_{\ell^2}.
\end{align}
In the supremum, the testing vectors may be confined to the same
finite coordinate space as $a$.

We will obtain
\eqref{eq:sr-dual-euclidean} by applying Stein's inequality in $E^*$,
which has the Stein property by the self-duality recalled in the
introduction.

\smallskip
{\em Step 1} --
We first need a bound for the functionals associated with disjoint
intervals. For integer intervals $H_1,\ldots,H_d$ listed from left
to right, we claim that
\begin{equation}\label{eq:sr-dual-interval}
 \Bigl\|\sum_{i=1}^d b_i h_{H_i}\Bigr\|_{E^*}
 \leq2AK\rho^*(b),
\end{equation}
where, as in \eqref{eq:defvHhH}, $H_i=[\lambda_i,\mu_i)$,  with
$1\leq\lambda_i<\mu_i$ and $\mu_i\leq\lambda_{i+1}$, and $h_{H_i} = \sum_{k=\lambda_i}^{\mu_i-1}e_k^*$.
For a finitely supported vector $x=\sum_k x_ke_k$, by \eqref{eq:defP-K} the truncation
$x':=P_{\mu_d-1}x$ to the first $\mu_d-1$ coordinates satisfies $\|x'\|_E\leq K\|x\|_E$.
It retains every coefficient occurring in an interval $H_i$,
so $h_{H_i}(x')=h_{H_i}(x)$.
Let
\[
 g_1:=\sum_{k=1}^{\lambda_1-1}x_k,\qquad c_i:=\sum_{k=\lambda_i}^{\mu_i-1}x_k,\qquad
 g_i:=\sum_{k=\mu_{i-1}}^{\lambda_i-1}x_k, \qquad 2\leq i\leq d.
\]
An empty sum is understood to be zero. In their natural order,
the gaps and the intervals $H_i$ partition the indices
$1,\ldots,\mu_d-1$.
The consecutive blocks are determined by the nondecreasing list
of cut points
\[
 0,\ \lambda_1-1,\ \mu_1-1,\ \lambda_2-1,\ \mu_2-1,\ldots,
 \lambda_d-1,\ \mu_d-1.
\]
Removing repeated cut points gives a strictly increasing list
$0=k_0<\cdots<k_m=\mu_d-1$. For each resulting nonempty
block, put
\[
 s_j:=\sum_{k=k_{j-1}+1}^{k_j}x_k,
 \qquad 1\leq j\leq m.
\]
Applying \eqref{eq:sr-compression} gives
\[
 \Bigl\|\sum_{j=1}^m s_j e_j\Bigr\|_E
 =\Bigl\|\sum_{j=1}^m
       \Bigl(\sum_{k=k_{j-1}+1}^{k_j}x_k\Bigr)e_j\Bigr\|_E
 \leq\Bigl\|\sum_{k=1}^{\mu_d-1}x_ke_k\Bigr\|_E
 =\|x'\|_E.
\]

Each block is either a gap, whose coefficient sum is $g_i$,
or an interval $H_i$, whose coefficient sum is $c_i$.
Thus the list $(s_1,\ldots,s_m)$ is
$(g_1,c_1,\ldots,g_d,c_d)$ with precisely those entries $g_i$
corresponding to empty gaps omitted. Each omitted entry is zero.
Reinserting these zeros preserves the norm by spreading
invariance. Consequently
\[
 \Bigl\|\sum_{i=1}^d(g_i e_{2i-1}+c_i e_{2i})\Bigr\|_E
 =\Bigl\|\sum_{j=1}^m s_j e_j\Bigr\|_E
 \leq\|x'\|_E\leq K\|x\|_E.
\]

Put
\[
 y:=\sum_{i=1}^d(g_i e_{2i-1}+c_i e_{2i}).
\]
The preceding estimate gives $\|y\|_E\leq K\|x\|_E$.
Let $M$ be the linear operator on
$\operatorname{span}\{e_1,\ldots,e_{2d}\}$ determined by
\[
 Me_{2i-1}=Me_{2i}=e_{2i-1},
 \qquad 1\leq i\leq d.
\]
For any $\xi=\sum_{k=1}^{2d}\xi_k e_k$, spreading invariance
and \eqref{eq:sr-compression}, with cut points
$0,2,\ldots,2d$, give
\[
 \|M\xi\|_E
 =\Bigl\|\sum_{i=1}^d(\xi_{2i-1}+\xi_{2i})e_i\Bigr\|_E
 \leq\|\xi\|_E.
\]
In particular,
\[
 My=\sum_{i=1}^d(g_i+c_i)e_{2i-1},
 \qquad \|My\|_E\leq\|y\|_E.
\]
Subtracting $y$ from $My$ cancels the gap coefficients $g_i$.
Consequently,
\begin{equation}\label{eq:sr-opposite-pairs}
 \begin{aligned}
 z:=My-y&=\sum_{i=1}^d c_i(e_{2i-1}-e_{2i}),\\
 \|z\|_E&\leq\|My\|_E+\|y\|_E
          \leq2\|y\|_E\leq2K\|x\|_E.
 \end{aligned}
\end{equation}

Suppressing the even coordinates of the coefficient vector of $z$
and using spreading invariance gives
$\rho(c)\leq\rho(c_1,-c_1,\ldots,c_d,-c_d)$.
Lemma~\ref{lem:sr-randomisation}, applied to the coefficient vector
\[
 \widetilde c:=(c_1,-c_1,\ldots,c_d,-c_d,0,\ldots)\in c_{00},
\]
gives
\[
 \rho(\widetilde c)
 \leq A|\widetilde c|_E
 =A\Bigl\|\sum_{i=1}^d c_i(e_{2i-1}-e_{2i})\Bigr\|_E
 =A\|z\|_E.
\]
Combining this with the preceding suppression estimate and
\eqref{eq:sr-opposite-pairs} gives
\[
 \rho(c)\leq\rho(\widetilde c)
 \leq A\|z\|_E
 \leq2AK\|x\|_E.
\]
Since $c_i=h_{H_i}(x)$, the defining dual inequality
\eqref{eq:sr-dual-norm} now gives
\[
 \Bigl|\sum_{i=1}^d b_i h_{H_i}(x)\Bigr|
 =\Bigl|\sum_{i=1}^d b_ic_i\Bigr|
 \leq\rho^*(b)\rho(c)
 \leq2AK\rho^*(b)\|x\|_E.
\]
The finite linear combination $\sum_i b_i h_{H_i}$ is a continuous
functional on $E$. Density of $c_{00}$ extends this estimate to
every $x\in E$, and taking the supremum over its unit ball proves
\eqref{eq:sr-dual-interval}.

\smallskip{\em Step 2} --
For $a=(a_1,\ldots,a_n)$, we use the random intervals constructed
at the beginning of this subsection and put
\begin{equation}\label{eq:sr-dual-test-inputs}
 g_j(\omega):=a_jh_{H_j(\omega)}
 \qquad 1\leq j\leq n,\ \omega\in\Omega_n,
\end{equation}
where again each $h_{H_j}$ is defined in \eqref{eq:defvHhH}.
Each function takes only finitely many values in $E^*$, so it
is an admissible input in \eqref{eq:PX-finite-Stein}.
For fixed $\omega\in\Omega_n$, let $\pi$ list the indices of the
intervals $H_j(\omega)$ from left to right. The coefficients in
that order are
\[b=(r_{\pi(1)}a_{\pi(1)},\,\ldots\,,r_{\pi(n)}a_{\pi(n)}).\]
The first part of this list consists of the left-half choices
in chronological order. The second consists of the right-half
choices in reverse chronological order.

For an explicit estimate, let
\[
 J_{\rm L}=\{i_1<\cdots<i_p\},
 \qquad J_{\rm R}=\{j_1<\cdots<j_q\}
\]
be the left-half and right-half indices, respectively.
Then $p+q=n$ and
\[
 (\pi(1),\ldots,\pi(n))
 =(i_1,\ldots,i_p,j_q,\ldots,j_1).
\]
Write $0_m$ for a list of $m$ zeros, and let $a_J$ denote
$a$ with all coordinates outside $J$ set to zero.
Using successively sign invariance, the triangle inequality,
spreading invariance to delete zeros, reversal invariance,
spreading invariance to restore the original coordinate positions,
and finally \eqref{eq:sr-dual-suppression}, we obtain
\begin{equation}\label{eq:sr-dual-reordering}
 \begin{aligned}
 \rho^*((r_{\pi(i)}a_{\pi(i)})_{i=1}^n)
 &=\rho^*(a_{i_1},\ldots,a_{i_p},a_{j_q},\ldots,a_{j_1})\\
 &\leq\rho^*(a_{i_1},\ldots,a_{i_p},0_q)
       +\rho^*(0_p,a_{j_q},\ldots,a_{j_1})\\
 &=\rho^*(a_{i_1},\ldots,a_{i_p})
       +\rho^*(a_{j_q},\ldots,a_{j_1})\\
 &=\rho^*(a_{i_1},\ldots,a_{i_p})
       +\rho^*(a_{j_1},\ldots,a_{j_q})\\
 &=\rho^*(a_{J_{\rm L}})+\rho^*(a_{J_{\rm R}})
 \leq2\rho^*(a).
 \end{aligned}
\end{equation}
An empty coefficient list has norm zero, so the computation
also covers $p=0$ or $q=0$.

Applying \eqref{eq:sr-dual-interval} to the intervals in this
spatial order, and then using \eqref{eq:sr-dual-reordering}, gives
\begin{equation}\label{eq:sr-dual-input}
 \Bigl\|\sum_{j=1}^n r_jg_j(\omega)\Bigr\|_{E^*}
 \leq4AK\rho^*(a).
\end{equation}
This estimate holds for every $\omega$ and every choice of the signs.

For the inputs \eqref{eq:sr-dual-test-inputs} we take $k_j=j-1$.
By \eqref{eq:sr-averages},
\[
 \mathcal E_{j-1}^{(n)}g_j(\omega)
 =\frac{1}{2}a_jh_{I_{j-1}(\omega)}.
\]
For fixed $\omega$, let $k(\omega):=\lambda_{n-1}(\omega)$ be the
left endpoint of the smallest parent interval.
Nesting gives $k(\omega)\in I_{j-1}(\omega)$ for all $j$.
The definition \eqref{eq:defvHhH} of $h_H$ consequently gives
$h_{I_{j-1}(\omega)}(e_{k(\omega)})=1$ for all $j$.
Since $\|e_{k(\omega)}\|_E=1$, the definition of the dual norm implies
\[
 \Bigl\|\frac12\sum_{j=1}^n r_ja_j
                   h_{I_{j-1}(\omega)}\Bigr\|_{E^*}
 \geq
 \Bigl|\Bigl(\frac12\sum_{j=1}^n r_ja_j
                   h_{I_{j-1}(\omega)}\Bigr)(e_{k(\omega)})\Bigr|
 =\frac12\Bigl|\sum_{j=1}^n r_ja_j\Bigr|.
\]
After squaring and averaging, independence and
$\E_r(r_ir_j)=\mathbf1_{\{i=j\}}$ give
\begin{equation}\label{eq:sr-dual-output}
 \begin{aligned}
 \E_\omega\E_r
 \Bigl\|\sum_{j=1}^n r_j\mathcal E_{j-1}^{(n)}g_j(\omega)
       \Bigr\|_{E^*}^2
 \geq\frac14\E_r\Bigl|\sum_{j=1}^n r_ja_j\Bigr|^2
 =\frac14\sum_{j=1}^n|a_j|^2.
 \end{aligned}
\end{equation}
Taking square roots in \eqref{eq:sr-dual-output}, applying
\eqref{eq:PX-finite-Stein} in $E^*$, and integrating the pointwise
input bound \eqref{eq:sr-dual-input}, we obtain
\begin{align*}
 \frac12\|a\|_{\ell^2}
 &\leq
 \Bigl\|\sum_{j=1}^n r_j\mathcal E_{j-1}^{(n)}g_j
       \Bigr\|_{L^2(\Omega_n\times\Omega_r;E^*)}\\
 &\leq\mathsf S(E^*)
 \Bigl\|\sum_{j=1}^n r_jg_j
       \Bigr\|_{L^2(\Omega_n\times\Omega_r;E^*)}
 \leq4AK\mathsf S(E^*)\rho^*(a).
\end{align*}
This proves \eqref{eq:sr-dual-euclidean} with
$B:=8AK\mathsf S(E^*)$, independently of $n$.
Equations \eqref{eq:sr-coefficient-duality} and
\eqref{eq:sr-dual-euclidean}, as explained at the start of the proof,
now give the asserted estimate \eqref{eq:rho-Brhostar} for $\rho$.
\end{proof}

\subsection{Proof of Theorem~\ref{thm:Stein-superreflexive}}

\begin{proof}[Proof of Theorem~\ref{thm:Stein-superreflexive}]
Suppose that a Stein space exists which is not superreflexive, and let $E$
be given by Lemma~\ref{lem:sr-reduction}.
Lemmas~\ref{lem:sr-randomisation} and~\ref{lem:sr-euclidean},
together with \eqref{eq:sr-lower}, give constants $A,B<\infty$,
independent of the number of coefficients, such that
\begin{equation}\label{eq:sr-two-bounds}
 \|a\|_{\ell^2}\leq\rho(a)\leq A|a|_E,
 \qquad
 \rho(a)\leq B\|a\|_{\ell^2}.
\end{equation}
For $a=(a_1,\ldots,a_N)$ and $0=k_0<\cdots<k_d=N$, put
$c_j:=\sum_{i=k_{j-1}+1}^{k_j}a_i$.
Applying the first two inequalities in \eqref{eq:sr-two-bounds}
to the shorter vector $c=(c_1,\ldots,c_d)$ gives
$\|c\|_{\ell^2}\leq\rho(c)\leq A|c|_E$.
The block-compression estimate \eqref{eq:sr-compression} gives
$|c|_E\leq|a|_E$. Hence
\begin{equation}\label{eq:sr-variation}
 \Bigl(\sum_{j=1}^d
       \Bigl|\sum_{i=k_{j-1}+1}^{k_j}a_i\Bigr|^2\Bigr)^{1/2}
 \leq A|a|_E = A\Bigl\|\sum_{i=1}^N a_i e_i\Bigr\|_E.
\end{equation}

Let $w$ be a continuous real function on $[0,1]$, and fix a
partition $0=t_0<\cdots<t_N=1$.
Put $a_i:=w(t_i)-w(t_{i-1})$ for $1\leq i\leq N$.

For every choice of indices $0=k_0<\cdots<k_d=N$,
telescoping and \eqref{eq:sr-variation} give
\begin{align*}
 \Bigl(\sum_{j=1}^d
       |w(t_{k_j})-w(t_{k_{j-1}})|^2\Bigr)^{1/2}
 =
 \Bigl(\sum_{j=1}^d
       \Bigl|\sum_{i=k_{j-1}+1}^{k_j}a_i\Bigr|^2\Bigr)^{1/2}
 \leq A\Bigl\|\sum_{i=1}^N a_i e_i\Bigr\|_E.
\end{align*}
The right-hand side depends only on the increments over the
original partition $0=t_0<\cdots<t_N=1$, and is independent
of the chosen indices $k_j$ and their number $d$.
Taking the maximum over all subpartitions and substituting the definition of $a_i$ therefore gives
\[
 \max_{\substack{1\leq d\leq N\\0=k_0<\cdots<k_d=N}}
 \Bigl(\sum_{j=1}^d
       |w(t_{k_j})-w(t_{k_{j-1}})|^2\Bigr)^{1/2}
 \leq
 A\Bigl\|\sum_{i=1}^N
       \bigl(w(t_i)-w(t_{i-1})\bigr)e_i\Bigr\|_E.
\]

Define the {\em strong $2$-variation}
\begin{equation}\label{eq:sr-strong-variation}
 V_2(w):=\sup_{d\ge 1}\sup_{0=t_0<\cdots<t_d=1}
       \Bigl(\sum_{j=1}^d|w(t_j)-w(t_{j-1})|^2\Bigr)^{1/2}.
\end{equation}
Let $W$ be standard real Brownian motion on $[0,1]$.
We use the classical path property
\begin{equation}\label{eq:sr-Brownian}
 V_2(W)=\infty\quad\hbox{almost surely};
\end{equation}
see \cite{Taylor-Brownian-variation} and
\cite[Section~13.9]{FrizVictoir}.
We will derive a contradiction, and thus prove the theorem,
by showing that \eqref{eq:sr-two-bounds} instead gives a finite
second moment $\mathbb E V_2(W)^2$.

For $m\geq0$, we restrict \eqref{eq:sr-strong-variation}
to the grid $D_m:=\{i2^{-m}:0\leq i\leq2^m\}$ and put
\begin{equation}\label{eq:sr-dyadic-variation-definition}
 V_{2,m}(w):=
 \max_{\substack{1\leq d\leq2^m\\0=k_0<\cdots<k_d=2^m}}
 \Bigl(\sum_{j=1}^d
       |w(k_j2^{-m})-w(k_{j-1}2^{-m})|^2\Bigr)^{1/2}.
\end{equation}
Every partition in the maximum is a subpartition of $D_m$.
Writing
\[\Delta_i^mW:=W(i2^{-m})-W((i-1)2^{-m}),\]
we may therefore apply \eqref{eq:sr-variation} to the coefficient
vector $(\Delta_i^mW)_{i=1}^{2^m}$.
The right-hand side is the same for every subpartition, so
taking their maximum gives
\begin{equation}\label{eq:sr-dyadic-variation}
 V_{2,m}(W)
 \leq A\Bigl\|\sum_{i=1}^{2^m}\Delta_i^mW\,e_i\Bigr\|_E.
\end{equation}

For fixed $m$, the increments $\Delta_i^mW$ are independent centred
Gaussian variables of variance $2^{-m}$.
Symmetry of a nondegenerate centred Gaussian density shows that
its sign is equally likely to be $-1$ or $1$, independently of its
absolute value. Independence of the increments gives this property
simultaneously for all their signs.
Conditionally on the absolute values of the increments, their
signs are therefore a Rademacher sequence.
By the definition \eqref{eq:def-rhoa},
\begin{equation}\label{eq:sr-increment-randomisation}
 \E\Bigl(\Bigl\|\sum_i\Delta_i^mW\,e_i\Bigr\|_E^2
       \,\Big|\, (|\Delta_i^mW|)_i\Bigr)
 =\rho((|\Delta_i^mW|)_i)^2.
\end{equation}
Let $\mathscr H_m:=\sigma(|\Delta_i^mW|:1\leq i\leq2^m)$.
Using \eqref{eq:sr-dyadic-variation}, the tower property of
conditional expectation, \eqref{eq:sr-increment-randomisation},
and the second bound in \eqref{eq:sr-two-bounds}, we obtain
\begin{equation}\label{eq:sr-variation-moment}
 \begin{aligned}
 \E V_{2,m}(W)^2
 &\leq A^2\E
       \Bigl\|\sum_{i=1}^{2^m}\Delta_i^mW\,e_i\Bigr\|_E^2
 =A^2\E\,\E\Bigl(
       \Bigl\|\sum_{i=1}^{2^m}\Delta_i^mW\,e_i\Bigr\|_E^2
       \,\Big|\,\mathscr H_m\Bigr)\\
 &=A^2\E\,\rho((|\Delta_i^mW|)_{i=1}^{2^m})^2
 \leq A^2B^2\E
       \bigl\|(|\Delta_i^mW|)_{i=1}^{2^m}\bigr\|_{\ell^2}^2\\
 &=A^2B^2\E\sum_{i=1}^{2^m}|\Delta_i^mW|^2
 =A^2B^2\,2^m2^{-m}=A^2B^2.
 \end{aligned}
\end{equation}
Here each increment is centred and has variance $2^{-m}$,
so $\E|\Delta_i^mW|^2=2^{-m}$.

We next pass from the dyadic grids to all partitions.
By \eqref{eq:sr-dyadic-variation-definition}, $V_{2,m}(W)$ is
the maximum of finitely many measurable functions and hence is
measurable. Since $D_m\subseteq D_{m+1}$, every subpartition
allowed at level $m$ is also allowed at level $m+1$. Therefore
\[
 V_{2,m}(w)\leq V_{2,m+1}(w)\leq V_2(w).
\]
In particular, the limit of $V_{2,m}(w)$ exists. To identify it,
we approximate an arbitrary partition by dyadic ones.

Fix any finite partition $0=t_0<\cdots<t_d=1$ and put
$t_j^{(m)}:=2^{-m}\lfloor2^m t_j\rfloor$.
For sufficiently large $m$, the number $2^{-m}$ is smaller than
every gap $t_j-t_{j-1}$. The approximating points are then strictly
increasing, belong to $D_m$, and have endpoints $0$ and $1$.
They also converge to the corresponding $t_j$.
For continuous $w$ this gives
\begin{align*}
 \Bigl(\sum_{j=1}^d|w(t_j)-w(t_{j-1})|^2\Bigr)^{1/2}
 =\lim_{m\to\infty}
   \Bigl(\sum_{j=1}^d
       |w(t_j^{(m)})-w(t_{j-1}^{(m)})|^2\Bigr)^{1/2}
 \leq\lim_{m\to\infty}V_{2,m}(w).
\end{align*}
The right-hand side is independent of the partition.
Taking the supremum in \eqref{eq:sr-strong-variation} proves
$V_{2,m}(w)\uparrow V_2(w)$.
Since Brownian paths are continuous almost surely, this estimate applies to $w=W$ almost surely.
It also identifies $V_2(W)$ with an increasing limit of measurable
random variables.
By monotone convergence and \eqref{eq:sr-variation-moment},
\[
 \E V_2(W)^2=\lim_{m\to\infty}\E V_{2,m}(W)^2
 \leq A^2B^2<\infty.
\]
Thus $V_2(W)$ is finite almost surely, contradicting
\eqref{eq:sr-Brownian}. This completes the proof of Theorem \ref{thm:Stein-superreflexive}.
\end{proof}

\section{Qiu's superreflexive space fails the Stein property}

\label{sec:TH-Qiu}
Fix $1<p\ne q<\infty$, and define the complex spaces
\begin{equation*}
 Y_0:=\mathbb C,\qquad
 Y_{N+1}:=L^p(\D;L^q(\D;Y_N)), \qquad N\geq0.
\end{equation*}
Thus $Y_N$ is an iterated mixed-norm space on $\D^{2N}$,
of dimension $4^N$. Its coordinate indicators form a
$1$-unconditional basis: multiplying any selection of their
coefficients by signs leaves the mixed norm unchanged.
Identifying $y\in Y_N$ with the constant function
$(u,v)\mapsto y$ gives an isometric inclusion $Y_N\to Y_{N+1}$.
We denote the completion of the resulting increasing union by
\begin{equation*}
 Y(p,q):=\overline{\bigcup_{N\geq0}Y_N}.
\end{equation*}
Qiu proves that this space is superreflexive
\cite[Corollary~5.3]{Qiu}.

To relate his estimates to the Stein property, let $F$ be a Banach
space with a fixed $1$-unconditional basis $(u_i)_{i\ge 1}$.
Following \cite[Definition~3.1]{Qiu}, write $S_{\mathrm Q}(F)$
for the least constant $C$ such that
\begin{equation}\label{eq:TH-Qiu-constant}
 \Bigl\|\sum_{j=1}^d(\mathcal E_j\varphi_j)u_{i_j}
       \Bigr\|_{L^1(\Omega;F)}
 \leq C\Bigl\|\sum_{j=1}^d\varphi_ju_{i_j}
       \Bigr\|_{L^\infty(\Omega;F)}.
\end{equation}
The estimate is required for every probability space $\Omega$,
every finite filtration $(\mathscr A_j)_{j=1}^d$ with
$\mathcal E_j:=\E(\,\cdot\mid\mathscr A_j)$, every list of
distinct basis indices $i_1,\ldots,i_d$, and every choice of
bounded scalar functions $\varphi_1,\ldots,\varphi_d$.

Suppose that $F$ has the Stein property. The version on arbitrary
probability spaces recalled in the introduction gives a constant
$C_F$ for its $L^2$ Stein inequality with finite filtrations.
For independent external Rademacher signs,
$1$-unconditionality gives, at every $\omega\in\Omega$,
\[
 \Bigl\|\sum_jr_j\varphi_j(\omega)u_{i_j}\Bigr\|_F
 =\Bigl\|\sum_j\varphi_j(\omega)u_{i_j}\Bigr\|_F.
\]
The same identity holds for the coefficients
$(\mathcal E_j\varphi_j)(\omega)$.
The inclusions $L^\infty(\Omega;F)\subseteq L^2(\Omega;F)
\subseteq L^1(\Omega;F)$ and the Stein inequality therefore give
\begin{equation}\label{eq:TH-Qiu-Stein-comparison}
 \begin{aligned}
 \Bigl\|\sum_j(\mathcal E_j\varphi_j)u_{i_j}\Bigr\|_{L^1(\Omega;F)}
 &\leq\Bigl\|\sum_j(\mathcal E_j\varphi_j)u_{i_j}
            \Bigr\|_{L^2(\Omega;F)}\\
 &=\Bigl(\E_r\Bigl\|\sum_jr_j(\mathcal E_j\varphi_j)u_{i_j}
            \Bigr\|_{L^2(\Omega;F)}^2\Bigr)^{1/2}\\
 &\leq C_F\Bigl(\E_r\Bigl\|\sum_jr_j\varphi_ju_{i_j}
            \Bigr\|_{L^2(\Omega;F)}^2\Bigr)^{1/2}\\
 &=C_F\Bigl\|\sum_j\varphi_ju_{i_j}\Bigr\|_{L^2(\Omega;F)}
 \leq C_F\Bigl\|\sum_j\varphi_ju_{i_j}\Bigr\|_{L^\infty(\Omega;F)}.
 \end{aligned}
\end{equation}
Comparing \eqref{eq:TH-Qiu-Stein-comparison} with
\eqref{eq:TH-Qiu-constant} shows that $S_{\mathrm Q}(F)\leq C_F$.

On the other hand, \cite[Theorem~3.10]{Qiu} gives a constant $c(p,q)>1$ such that
\begin{equation*}
 S_{\mathrm Q}(Y_N)\geq c(p,q)^N,\qquad N\geq1.
\end{equation*}
where the constants are computed with respect to the coordinate
bases. If $Y(p,q)$ had the Stein property with a full $L^2$
Stein constant $C_Y$, its restriction to the isometric subspace
$Y_N$ would give the same bound for all $Y_N$-valued tests.
Applying \eqref{eq:TH-Qiu-Stein-comparison} in each $Y_N$ would
then give
\[
 c(p,q)^N\leq S_{\mathrm Q}(Y_N)\leq C_Y,
 \qquad N\geq1,
\]
which is impossible. Thus $Y(p,q)$ is superreflexive and fails
the Stein property.

\section[A twisted Hilbert space without the Stein property]{A twisted space without the Stein property}
\label{sec:twisted-Hilbert}

In this section we construct an example of a {\em twisted Hilbert space}
without the Stein property, and use it to conclude that the Stein property is not the three-space property.

The construction can be outlined as follows.
For each $N\geq1$, we construct a finite-dimensional complex
Hilbert space $H_N$ and a Banach norm on the vector space
$Z_N:=H_N\times H_N$. Using the maps
\[
 j_Nu:=(u,0),\qquad \pi_N(u,x):=x,
\]
this gives a short exact sequence
\[
 0\longrightarrow H_N\xrightarrow{\ j_N\ }Z_N
   \xrightarrow{\ \pi_N\ }H_N\longrightarrow0.
\]
Recall that an {\em exact sequence} means that every arrow represents a bounded operator whose range is closed and equals the kernel of the operator represented by the next arrow. A {\em short exact sequence} is an exact of five spaces where the initial and terminal space are zero; in such a sequence, the two outer mappings are the trivial mappings, the second is injective, and the third surjective. In the present example, this means that both the closed subspace $j_N(H_N)$
and the quotient $Z_N/j_N(H_N)$ will be isomorphic to $H_N$.
The block norms will satisfy
\[
 \|j_Nu\|_{Z_N}=\tfrac14\|u\|_{H_N},
 \qquad \|\pi_Nz\|_{H_N}\leq\|z\|_{Z_N},
\]
and every $x\in H_N$ will have a preimage $z_x\in Z_N$ with
$\|z_x\|_{Z_N}=\|x\|_{H_N}$.
Consequently, the direct sums
\[
 H:=\Bigl(\bigoplus_{N\geq1}H_N\Bigr)_{\ell^2},
 \qquad
 Z:=\Bigl(\bigoplus_{N\geq1}Z_N\Bigr)_{\ell^2}
\]
will fit into an exact sequence
\[
 \begin{gathered}
 0\longrightarrow H\xrightarrow{\ j\ }Z
   \xrightarrow{\ \pi\ }H\longrightarrow0,\\
 j((u_N)_N):=((u_N,0))_N,
 \qquad \pi((u_N,x_N)_N):=(x_N)_N.
 \end{gathered}
\]
Since $H$ is a Hilbert space, this will imply superreflexivity
of $Z$ by the three-space theorem recalled below.

To prove that $Z$ fails the Stein property, we will construct
explicit Stein tests in $Z_N$.
Put $I:=\{-1,1\}\times\{0,1\}$. For each $N$, we will construct
a finite probability space $W_N$, an increasing list of
$\sigma$-algebras $(\mathcal G_\alpha^{(N)})_{\alpha\in I^N}$,
and functions $f_\alpha:W_N\to Z_N$ such that
\begin{align*}
 \Bigl\|\sum_{\alpha\in I^N}r_\alpha f_\alpha
       \Bigr\|_{L^2(W_N\times\Omega_r;Z_N)}
 &=1,\\
 \Bigl\|\sum_{\alpha\in I^N}r_\alpha
       \E[f_\alpha\mid\mathcal G_\alpha^{(N)}]
       \Bigr\|_{L^2(W_N\times\Omega_r;Z_N)}
 &\geq b(1+\kappa\sqrt N),
\end{align*}
where $b,\kappa>0$ are independent of $N$.
The variables $r_\alpha$ are independent Rademacher signs,
independent of the variable in $W_N$.
Thus the ratio of the output norm to the input norm tends
to infinity. The isometric inclusions $Z_N\to Z$ carry these
tests into $Z$, proving that $Z$ fails the Stein property.

Throughout this section all Banach and Hilbert spaces are complex,
and the finite probability spaces used in our construction
are equipped with uniform probability measure.

\subsection{Statement of the results}
For a Banach space $X$ and a fixed number of summands $n\geq1$, put
\begin{equation}\label{eq:TH-finite-Stein-constants}
 \mathfrak s_n(X):=
 \inf\bigl\{C\geq0:\,
    \text{\eqref{eq:PX-finite-Stein} holds with $C_X=C$}\bigr\}.
\end{equation}
The estimate \eqref{eq:PX-finite-Stein} is required to hold uniformly with respect to $d\ge 1$ and the choices of
$k_1,\ldots,k_n\in\{0,\ldots,d\}$ and inputs $f_1,\ldots,f_n \in L^2(\{-1,1\}^d;X)$.

We observe the trivial bound $$\mathfrak s_n(X)\leq n.$$
Indeed, if
$F=\sum_{i=1}^n r_i f_i$, then $f_j=\E_r(r_jF)$.
The triangle inequality and the Cauchy--Schwarz inequality give
$\|f_j\|_{L^2}\leq\|F\|_{L^2(\Omega_r\times\D^d;X)}$.
Contractivity of conditional expectation and the triangle inequality
therefore give the asserted bound.

By adding zero terms we see that $\mathfrak s_1(X) \le \mathfrak s_2(X) \le \ldots$, and testing on a
single nonzero constant term gives the lower bound $1\le \mathfrak s_1(X)$ when
$X\ne\{0\}$. In summary, when
$X\ne\{0\}$, for all $n\ge 1$ we have
\begin{equation}\label{eq:TH-finite-Stein-basic}
 1\leq\mathfrak s_n(X)\leq n,\qquad
 \mathfrak s_n(X)\leq\mathfrak s_{n+1}(X).
\end{equation}
By \eqref{eq:TH-finite-Stein-constants} and the definition
\eqref{eq:PX-finite-Stein},
\begin{equation}\label{eq:TH-finite-Stein-criterion}
 X\text{ has the Stein property}
 \quad\Longleftrightarrow\quad
 \sup_{n\geq1}\mathfrak s_n(X)<\infty.
\end{equation}

\begin{definition}
A \emph{twisted Hilbert space} is a Banach space $X$ containing a
closed subspace $H$ isomorphic to a Hilbert space such that the quotient space $X/H$ is
isomorphic to a Hilbert space.
\end{definition}

In this situation $$0 \longrightarrow H \longrightarrow X \longrightarrow X/H \longrightarrow 0$$ is a short exact sequence. When admits a direct sum decomposition of the form $X = H \oplus E$, with $E$ isomorphic to $X/H$, the short exact sequence is said to be {\em split}. In this case $X$ is isomorphic to a Hilbert space. From our perspective this is the uninteresting situation. In the non-split case it may happen that $X$ is a Banach space which is not isomorphic to a Hilbert space; this will be the case in our example.

We will use the three-space property for superreflexivity:
if $Y$ is a closed subspace of a Banach space $X$ and both
$Y$ and $X/Y$ are superreflexive, then $X$ is superreflexive
\cite{EnfloLindenstraussPisier1975}.
In particular, every twisted Hilbert space is superreflexive,
because Hilbert spaces are superreflexive and superreflexivity
is invariant under isomorphisms.

The main result of this section reads as follows.

\begin{theorem}[A twisted Hilbert space without the Stein property]
\label{thm:TH-main}
There are a Hilbert space $H$, a superreflexive Banach space $Z$,
and an exact sequence of bounded linear maps
\[
 0\longrightarrow H\longrightarrow Z\longrightarrow H
 \longrightarrow0
\]
such that, for a numerical constant $c>0$,
\begin{equation}\label{eq:TH-main-bounds}
 \mathfrak s_n(Z)\geq c\sqrt{\log(n+1)}
 \qquad(n\geq2).
\end{equation}
In particular, $Z$ is a twisted Hilbert space, hence superreflexive, and does not have
the Stein property.
\end{theorem}

A Lean~4 formalisation of Theorem~\ref{thm:TH-main} is included
in the ancillary files. It verifies the construction of the twisted
Hilbert space and the quantitative lower bound for its finite Stein
constants. Superreflexivity is expressed through the equivalent
condition that the space admits an equivalent uniformly nonsquare
norm.

\medskip
Kalton \cite[Theorem~9.11]{Kalton-differentials} showed that an
extension of a Hilbert space by a Hilbert space can fail to be UMD.
Corollary~\ref{cor:TH-three-space} strengthens this conclusion:
even Bourgain's vector-valued Stein inequality, which holds in
every UMD space, can fail in such an extension. To prove this for
the space $Z$ of Theorem~\ref{thm:TH-main}, we construct explicit
tests for Stein's inequality in its finite-dimensional summands $Z_N$.

Theorem~\ref{thm:TH-main} has the following consequence, whose proof is included at the end of the section.

\begin{corollary}\label{cor:TH-three-space}
The Stein property is not a three-space property, even when
both the closed subspace and the quotient are Hilbert spaces.
\end{corollary}

The space in the theorem is of the form $Z = (\bigoplus_{N\ge 1} Z_N)_{\ell^2}$, where the spaces $Z_N$ have the properties described at the beginning of the section.
The construction of these spaces starts with a test on four coordinates. We iterate
this test using a filtration modelled on the one in Qiu's proof of
multiplicativity \cite[Theorem~3.5]{Qiu}. At each step, the
coordinatewise conditional expectations produce a new term orthogonal
to the contribution from the preceding steps. Its squared Hilbert
norm is $(\log 9)^2/5$. Thus, after $N$ steps, the accumulated term
has norm $\sqrt{N}\log 9/\sqrt{5}$, while the randomised input
still has norm one.
The filtration that combines these levels follows the arrangement
in Qiu's multiplicativity argument \cite[Theorem~3.5]{Qiu}.

\subsection{Mixed norms and their derivatives}

Let $I=\{-1,1\}\times\{0,1\}$ with product uniform probability.
For $N\geq1$, functions on $I^N$ have variables ordered as
$(u_1,v_1,\ldots,u_N,v_N)$.
For $-1/4\leq t\leq1/4$, we set
\begin{equation}\label{eq:TH-mixed-exponents}
 \frac1{p_j(t)}:=\frac12+(-1)^{j+1}t
 \qquad(1\leq j\leq2N)
\end{equation}
and define the mixed-norm spaces
\begin{equation}\label{eq:TH-mixed-norms}
 E_N(t):=
 L^{p_1(t)}_{u_1}\bigl(L^{p_2(t)}_{v_1}
   (\cdots L^{p_{2N-1}(t)}_{u_N}
                (L^{p_{2N}(t)}_{v_N})\cdots)\bigr).
\end{equation}
Each two-point variable is integrated with probability $1/2$
at each of its values, starting with the innermost variable $v_N$.
For example, at $N=1$,
\[
 \|x\|_{E_1(t)}
 =\Bigl(\frac12\sum_{u\in \{-1,1\}}
       \Bigl(\frac12\sum_{v\in \{0,1\}}|x(u,v)|^{p_2(t)}
       \Bigr)^{p_1(t)/p_2(t)}\Bigr)^{1/p_1(t)}.
\]
At $t=0$ all the exponents in \eqref{eq:TH-mixed-exponents} equal two.
Repeatedly summing squares in \eqref{eq:TH-mixed-norms}
then gives the Hilbert space $H_N:=L^2(I^N)$, with
\begin{equation}\label{eq:TH-Hilbert-norm}
 \|x\|_{H_N}^2=4^{-N}\sum_{\alpha\in I^N}|x_\alpha|^2.
\end{equation}
We use the inner product
$\langle f,g\rangle_{H_N}:=4^{-N}\sum_{\alpha\in I^N}
 f_\alpha\overline{g_\alpha}$, linear in its first entry.
Since $1/p_j(t)+1/p_j(-t)=1$ for every $j$, successive applications
of H\"older's inequality, starting with $v_N$, give
\begin{equation}\label{eq:TH-mixed-holder}
 |\langle f,g\rangle_{H_N}|
 \leq\|f\|_{E_N(t)}\|g\|_{E_N(-t)}
 \qquad(|t|\leq1/4).
\end{equation}

We temporarily write the $2N$ variables as $\xi_1,\ldots,\xi_{2N}$.
For $x\in H_N$, let $a_j(x)$ be the average of $|x|^2$ over
$\xi_j,\ldots,\xi_{2N}$, keeping the first $j-1$ variables fixed.
If $\E_{\xi_j}$ denotes averaging only in the $j$-th variable,
this means that
\begin{equation}\label{eq:TH-square-averages}
 \begin{aligned}
 a_j(x)&:=\E_{\xi_j}\cdots\E_{\xi_{2N}}|x|^2,
       && 1\leq j\leq2N;\\
 a_{2N+1}(x)&:=|x|^2.
 \end{aligned}
\end{equation}
In particular, $a_1(x)=\|x\|_{H_N}^2$ by
\eqref{eq:TH-Hilbert-norm}. Equation \eqref{eq:TH-square-averages} gives
\begin{equation}\label{eq:TH-average-recursion}
 \E_{\xi_j}a_{j+1}(x)=a_j(x).
\end{equation}
At any coordinate $\alpha$ where $x_\alpha\ne0$, every
$a_j(x)$ evaluated along that coordinate is strictly positive.
We may therefore define
\begin{equation}\label{eq:TH-Omega}
 \Omega_N(x)_\alpha
 :=x_\alpha\sum_{j=1}^{2N}(-1)^{j+1}
              \log\frac{a_{j+1}(x)(\alpha)}{a_j(x)(\alpha)},
 \qquad x_\alpha\ne0.
\end{equation}
At a coordinate with $x_\alpha=0$ we set $\Omega_N(x)_\alpha=0$.
In particular $\Omega_N(0)=0$.
We call this nonlinear map the \emph{centralizer} of the scale.
For a scalar $\zeta\ne0$, all averages for $\zeta x$ are the
corresponding averages for $x$ multiplied by $|\zeta|^2$.
The factors cancel in the ratios, which proves
$\Omega_N(\zeta x)=\zeta\Omega_N(x)$.
The same identity holds at $\zeta=0$ by definition.

To estimate $\Omega_N$, we realise it as the derivative at zero
of an explicit analytic function. For background on mixed-norm
interpolation and its derivatives, see
\cite[Section~5.1]{BerghLofstrom} and
\cite{RochbergWeiss,Kalton-differentials}.
For $x\ne0$, put
\begin{equation}\label{eq:TH-extremal-ratios}
 b_j:=a_j(x)^{1/2},\qquad
 \tau_j:=
 \begin{cases}
 b_{j+1}/b_j,&b_j>0,\\
 0,&b_j=0.
 \end{cases}
\end{equation}
The first definition applies for $1\leq j\leq2N+1$ and the
second for $1\leq j\leq2N$.
Along a coordinate where $x_\alpha\ne0$, all these ratios are
strictly positive. At such a coordinate also put
$\operatorname{sgn}(x_\alpha):=x_\alpha/|x_\alpha|$.
We define the analytic function coordinatewise by
\begin{equation*}
 F_x(z)_\alpha:=
 \operatorname{sgn}(x_\alpha)b_1
       \prod_{j=1}^{2N}\tau_j(\alpha)^{1+2(-1)^{j+1}z},
 \qquad |\operatorname{Re}z|\leq\frac14,
\end{equation*}
at nonzero coordinates, and set it identically zero at the
remaining coordinates. We also set $F_0=0$.
Powers of positive numbers in this formula use their real logarithms.
Their absolute values are bounded on the closed strip, so $F_x$
is a bounded continuous vector-valued function, analytic inside.
The ratios in \eqref{eq:TH-extremal-ratios} telescope, giving
$\prod_j\tau_j=b_{2N+1}/b_1=|x|/b_1$.
Consequently $F_x(0)=x$, and differentiation gives
\begin{equation}\label{eq:TH-extremal-derivative}
 F_x'(0)
 =x\sum_{j=1}^{2N}2(-1)^{j+1}\log\tau_j
 =\Omega_N(x).
\end{equation}
Both identities also hold at the zero coordinates.

The function is norming at every vertical line of the strip.
To check this, put $\gamma_j(t):=1+2(-1)^{j+1}t$.
By \eqref{eq:TH-mixed-exponents}, $\gamma_j(t)p_j(t)=2$.
Imaginary powers of $\tau_j$ have modulus one, so
$|F_x(t+is)|=b_1\prod_j\tau_j^{\gamma_j(t)}$ at nonzero
coordinates.
When $a_j(x)>0$, \eqref{eq:TH-extremal-ratios} and
\eqref{eq:TH-average-recursion} give
\[
 \E_{\xi_j}\tau_j^2
 =\frac{\E_{\xi_j}a_{j+1}(x)}{a_j(x)}=1.
\]
In taking the innermost $L^{p_{2N}(t)}$-norm, all preceding
factors are constant in $\xi_{2N}$. The last factor therefore
contributes
$(\E_{\xi_{2N}}\tau_{2N}^2)^{1/p_{2N}(t)}=1$.
After that integration the same calculation applies to
$\tau_{2N-1}$, and continuing outwards removes all factors.
If an average is zero, every entry below that choice of the
preceding variables is zero and contributes zero to the norm.
When the average is positive but a later one vanishes, the
corresponding zero ratio is raised to a positive power because
$1/2\leq\gamma_j(t)\leq3/2$.
Thus the integrations also apply to vectors with zero coordinates.
Only $b_1=\|x\|_{H_N}$ remains, proving
\begin{equation}\label{eq:TH-norming}
 \|F_x(t+is)\|_{E_N(t)}=\|x\|_{H_N},
 \qquad |t|\leq1/4,\ s\in\mathbb R.
\end{equation}

For $x,y\in H_N$ consider the scalar function
\[
 \Phi_{x,y}(z):=\langle F_x(z),F_y(-\overline z)\rangle_{H_N}.
\]
It is analytic in the open strip, since conjugation in the second
entry makes each $\overline{F_y(-\overline z)_\alpha}$ analytic
in $z$. If $\operatorname{Re}z=t$, then
$\operatorname{Re}(-\overline z)=-t$. Thus \eqref{eq:TH-mixed-holder}
and \eqref{eq:TH-norming} give
\[|\Phi_{x,y}(z)|\leq\|x\|_{H_N}\|y\|_{H_N}\]
throughout the strip.
By \eqref{eq:TH-extremal-derivative}, differentiation at zero gives
\[
\Phi_{x,y}'(0)=\langle\Omega_N(x),y\rangle_{H_N}-\langle x,\Omega_N(y)\rangle_{H_N}.\]
Applying Cauchy's estimate on $|z|=r$ and letting $r\uparrow1/4$, we obtain
\begin{equation}\label{eq:TH-duality-bound}
 \big|\langle\Omega_N(x),y\rangle_{H_N}
       -\langle x,\Omega_N(y)\rangle_{H_N}\big|
 \leq 4\|x\|_{H_N}\|y\|_{H_N}.
\end{equation}

We now define the block norm on $Z_N:=H_N\times H_N$ by
\begin{equation}\label{eq:TH-block-norm}
 \|(u,x)\|_{Z_N}:=
 \max\left\{\|x\|_{H_N},\
 \frac14\sup_{\|y\|_{H_N}\leq1}
 \big|\langle u,y\rangle_{H_N}-\langle x,\Omega_N(y)\rangle_{H_N}\big|\right\}.
\end{equation}
For each fixed $y$ the map
$(u,x)\mapsto\langle u,y\rangle_{H_N}-\langle x,\Omega_N(y)\rangle_{H_N}$
is complex linear. By adding and subtracting
$\langle\Omega_N(x),y\rangle_{H_N}$ and applying
\eqref{eq:TH-duality-bound}, we see that its modulus is bounded by
\[(\|u-\Omega_N(x)\|_{H_N}+4\|x\|_{H_N})\|y\|_{H_N}.\]
Thus the supremum is finite and \eqref{eq:TH-block-norm} defines a seminorm. If it vanishes,
then $x=0$ and $\langle u,y\rangle_{H_N}=0$ for every $y$ in
the unit ball, whence $u=0$. Thus $Z_N$ is a normed space,
and it is complete by finite dimensionality. Its definition gives
\begin{equation}\label{eq:TH-quotient-lower}
 \|x\|_{H_N}\leq\|(u,x)\|_{Z_N}.
\end{equation}
On the kernel of $\pi_N$ and on the graph of $\Omega_N$, the norm takes the values
\begin{equation}\label{eq:TH-block-identities}
 \|(u,0)\|_{Z_N}=\tfrac14\|u\|_{H_N},
 \qquad
 \|(\Omega_N(x),x)\|_{Z_N}=\|x\|_{H_N}.
\end{equation}
The first follows from Hilbert space duality. For the second,
\eqref{eq:TH-duality-bound} bounds the supremum in
\eqref{eq:TH-block-norm} by $4\|x\|_{H_N}$, and the first
term in the maximum is $\|x\|_{H_N}$.

To work with this norm concretely, put
\begin{equation}\label{eq:TH-qnorm}
 q_N(u,x):=\|u-\Omega_N(x)\|_{H_N}+\|x\|_{H_N}.
\end{equation}
Writing
\[(u,x)=(u-\Omega_N(x),0)+(\Omega_N(x),x)\] and applying
\eqref{eq:TH-block-identities} gives
\[\|(u,x)\|_{Z_N}\leq \tfrac14\|u-\Omega_N(x)\|_{H_N}+\|x\|_{H_N} \leq q_N(u,x).\]
Conversely,
\[
 \tfrac14\|u-\Omega_N(x)\|_{H_N}
 =\|(u,x)-(\Omega_N(x),x)\|_{Z_N}
 \leq\|(u,x)\|_{Z_N}+\|x\|_{H_N}
 \leq 2\|(u,x)\|_{Z_N}
\]
by \eqref{eq:TH-quotient-lower}. Multiplying by four and adding $\|x\|_{H_N}\leq\|(u,x)\|_{Z_N}$ gives
\begin{equation}\label{eq:TH-uniform-equivalence}
 bq_N(u,x)\leq\|(u,x)\|_{Z_N}\leq q_N(u,x), \quad b:=1/9.
\end{equation}
In particular, $q_N$ is a quasinorm equivalent to the Banach norm.
Its homogeneity follows from homogeneity of $\Omega_N$, and the
triangle inequality for the Banach norm gives
$q_N(z+z')\leq b^{-1}(q_N(z)+q_N(z'))$.
All comparison constants are independent of $N$.

\subsection{The exact sequence and superreflexivity}
\label{subsec:TH-exact-sequence}

We define the Banach and Hilbert direct sums
\begin{equation}\label{eq:TH-sums}
 Z:=\Bigl(\bigoplus_{N\geq1}Z_N\Bigr)_{\ell^2},
 \qquad
 H:=\Bigl(\bigoplus_{N\geq1}H_N\Bigr)_{\ell^2},
\end{equation}
with norms
\[
 \|(z_N)_N\|_Z^2:=\sum_N\|z_N\|_{Z_N}^2,
 \qquad
 \|(x_N)_N\|_H^2:=\sum_N\|x_N\|_{H_N}^2.
\]
With the direct sums $H,Z$ from \eqref{eq:TH-sums}, define
\begin{equation}\label{eq:TH-injection}
 j:H\to Z,\qquad j((u_N)_N):=((u_N,0))_{N\geq1}.
\end{equation}
Squaring the first identity in \eqref{eq:TH-block-identities}
and summing over $N$ gives
\begin{equation}\label{eq:TH-kernel-embedding}
 \|j(u)\|_Z=\tfrac14\|u\|_H.
\end{equation}
Thus $j$ is a bounded linear injection with closed image, since $H$ is complete.

For $z=((u_N,x_N))_N\in Z$, put
\begin{equation}\label{eq:TH-quotient-map}
 \pi z:=(x_N)_N.
\end{equation}
The estimate \eqref{eq:TH-quotient-lower} gives
\[
 \sum_N\|x_N\|_{H_N}^2
 \leq\sum_N\|(u_N,x_N)\|_{Z_N}^2=\|z\|_Z^2.
\]
It proves at once that $\pi z\in H$ and that the linear map
$\pi:Z\to H$ is contractive.
Definitions \eqref{eq:TH-injection} and \eqref{eq:TH-quotient-map}
give $\pi j=0$.
Conversely, if $\pi z=0$, every $x_N$ is zero.
By \eqref{eq:TH-block-identities},
\[\sum_N\|u_N\|_{H_N}^2=16\|z\|_Z^2<\infty.\]
Hence $u=(u_N)_N$ belongs to $H$ and $z=j(u)$.
We have proved $\ker\pi=j(H)$.

For surjectivity, let $x=(x_N)\in H$.
The second identity in \eqref{eq:TH-block-identities} gives
\[
 \|(\Omega_N(x_N),x_N)\|_{Z_N}=\|x_N\|_{H_N}.
\]
Thus \[z_x:=((\Omega_N(x_N),x_N))_N\] belongs to $Z$,
$\pi z_x=x$, and $\|z_x\|_Z=\|x\|_H$. Since $\pi$ is
contractive, every preimage of $x$ has norm at least $\|x\|_H$.
The quotient norm associated with \eqref{eq:TH-quotient-map} therefore equals
\begin{equation}\label{eq:TH-quotient-isometry}
 \inf_{\pi z=x}\|z\|_Z=\|x\|_H.
\end{equation}
Thus $\pi$ induces an isometric isomorphism $Z/j(H)\to H$,
and
\[
 0\longrightarrow H\xrightarrow{\ j\ }Z
   \xrightarrow{\ \pi\ }H\longrightarrow0
\]
is the required exact sequence.

We now apply the three-space theorem stated above.
The subspace $j(H)$ is closed by \eqref{eq:TH-kernel-embedding}
and the argument following it, and is isomorphic to the Hilbert
space $H$. It is therefore superreflexive.
The quotient $Z/j(H)$ is isometric to $H$ by
\eqref{eq:TH-quotient-isometry}, so it too is superreflexive.
Consequently, the result of \cite{EnfloLindenstraussPisier1975}
shows that $Z$ is superreflexive.

\subsection{A four-coordinate seed}

Having proved that $Z$ is a superreflexive twisted Hilbert space,
we now construct the tests for the Stein lower bound
\eqref{eq:TH-main-bounds}.
The four coordinates of $H_1=L^2(I)$ will be assigned conditional
expectations in the order
\begin{equation}\label{eq:TH-coordinate-order}
 (-1,0)\prec(1,0)\prec(-1,1)\prec(1,1).
\end{equation}
This order specifies when information is revealed in the filtration.
The mixed norm $E_1(t)=L^{p_1(t)}_u(L^{p_2(t)}_v)$ continues to
integrate in $v$ first and then in $u$, as defined in
\eqref{eq:TH-mixed-norms}.

Let $w$ be a uniform sign and set
\begin{equation}\label{eq:TH-seed}
 a:=\sqrt{2/5},\qquad d:=4/\sqrt5,\qquad
 \theta_{u,0}(w):=a,\qquad
 \theta_{u,1}(w):=d\,\mathbf1_{\{u=w\}}.
\end{equation}
On the space of the sign $w$, put
\begin{equation}\label{eq:TH-seed-filtration}
 \mathcal A_{(u,0)}:=\{\emptyset,\D\},\qquad
 \mathcal A_{(u,1)}:=\sigma(w).
\end{equation}
The first two entries of the ordered list therefore use trivial
information and the last two use full information.
We define the linear operator $T_1$ on $H_1$-valued functions of $w$
by
\begin{equation}\label{eq:TH-seed-operator}
 (T_1f)_i:=\E[f_i\mid\mathcal A_i].
\end{equation}
At $v=0$ the coordinate $\theta_{u,0}$ is constant, so averaging
fixes it. At $v=1$ the conditional expectation uses the full
$\sigma$-algebra, so it fixes $\theta_{u,1}$ as well.
Thus $T_1\theta=\theta$.

We compute the part of $\Omega_1(\theta)$ that $T_1$ changes.
By \eqref{eq:TH-Hilbert-norm}, the Hilbert norm uses weight
$1/4$ on each coordinate, and
\[
 a_1(\theta)=\|\theta(w)\|_{H_1}^2
 =\frac{2a^2+d^2}{4}=1.
\]
The average \eqref{eq:TH-square-averages} in the inner variable is
\[
 a_2(\theta)(u)
 =\frac12\bigl(a^2+d^2\mathbf1_{\{u=w\}}\bigr)
 =\begin{cases}9/5,&u=w,\\1/5,&u=-w.\end{cases}
\]
At $v=0$ the final average is $a_3(\theta)(u,0)=a^2$.
The two logarithms in \eqref{eq:TH-Omega} consequently give
\[
 (\Omega_1(\theta))_{u,0}
 =a\Bigl(\log\frac{a_2(\theta)(u)}1
               -\log\frac{a^2}{a_2(\theta)(u)}\Bigr)
 =a\bigl(2\log a_2(\theta)(u)-\log a^2\bigr).
\]
The two values of $a_2(\theta)(u)$ can be written as
$a_2(\theta)(u)=(3/5)9^{uw/2}$.
Substitution shows that
\[
 (\Omega_1(\theta))_{u,0}
 =a\bigl(2\log(3/5)-\log a^2\bigr)+a\,uw\log9.
\]
The first term is independent of $w$, and the second has mean zero.
Thus taking its expectation and subtracting the original expression
leaves $-a\,uw\log9$.
At $v=1$, $T_1$ is the identity, so the corresponding difference
is zero, including at a zero coordinate of $\theta$.
With
\begin{equation}\label{eq:TH-seed-residual-definition}
 c:=T_1\Omega_1(\theta)-\Omega_1(\theta),
\end{equation}
where $T_1$ is given by \eqref{eq:TH-seed-operator}, we obtain
\begin{equation*}
 c_{u,0}(w)=-a\,uw\log9,\qquad c_{u,1}(w)=0.
\end{equation*}

We put
\begin{equation}\label{eq:TH-seed-kappa}
 \kappa:=\frac{\log9}{\sqrt5}.
\end{equation}
The inner product is linear in its first entry, so in general
$\langle x,y\rangle_{H_1}=\frac14\sum_i x_i\overline{y_i}$.
All entries of $\theta$ and $c$ are real.
Only the two entries with $v=0$ contribute to their inner product,
and their contributions cancel because $u$ takes the values
$-1$ and $1$. Explicitly,
\[
 \langle\theta(w),c(w)\rangle_{H_1}
 =-\frac{a^2w\log9}{4}\sum_{u=\pm1}u=0.
\]
Similarly,
\[
 \|c(w)\|_{H_1}^2
 =\frac14\sum_{u=\pm1}a^2u^2w^2(\log9)^2
 =\frac{(\log9)^2}{5}=\kappa^2.
\]
We have therefore proved, for both choices of $w$,
\begin{equation}\label{eq:TH-seed-identities}
 \|\theta(w)\|_{H_1}^2=1,\qquad
 \langle\theta(w),c(w)\rangle_{H_1}=0,\qquad
 \|c(w)\|_{H_1}^2=\kappa^2.
\end{equation}

\subsection{The branching filtration}

We iterate the seed using four independent copies at each level.
Let $H_0=\mathbb C$, let $W_0$ be a singleton with probability one,
and set $h_0=1$.
For $N\geq1$, let
\begin{equation}\label{eq:TH-branch-spaces}
 W_N:=\D\times\prod_{i\in I}W_{N-1}^{(i)}
\end{equation}
with product probability, where the superscript labels independent
copies of $W_{N-1}$.
A point is written $w=(w_0,(w_i)_{i\in I})$.
The coordinate $w_0$ is the new sign, while each $w_i$ contains
the signs used in the copy indexed by $i$.
We define $h_N:W_N\to H_N$ recursively by
\begin{equation}\label{eq:TH-h-recursion}
 (h_N(w))_{i,\alpha}
 :=\theta_i(w_0)(h_{N-1}(w_i))_\alpha,
 \qquad i\in I,\quad\alpha\in I^{N-1}.
\end{equation}
The product measure on $I^N=I\times I^{N-1}$ gives
\[
 \|h_N(w)\|_{H_N}^2
 =\frac14\sum_{i\in I}|\theta_i(w_0)|^2
                   \|h_{N-1}(w_i)\|_{H_{N-1}}^2.
\]
Starting with $\|h_0\|=1$, induction and
\eqref{eq:TH-seed-identities} show that the right-hand side is one
at every $w$. Thus
\begin{equation}\label{eq:TH-h-normalisation}
 \|h_N(w)\|_{H_N}=1\qquad(w\in W_N).
\end{equation}

The coordinates of $H_N$ are indexed by $I^N$.
We order this set lexicographically using
\eqref{eq:TH-coordinate-order} on each factor.
For an index $(i,\alpha)$, the information revealed consists of
the $\sigma$-algebra $\mathcal A_i$ of the root sign, all signs in
earlier copies $W_{N-1}^{(j)}$ with $j\prec i$, and the information
at stage $\alpha$ in the current copy.
If $\operatorname{pr}_0,\operatorname{pr}_j$ are the coordinate
projections from $W_N$ and $\mathscr P(W_{N-1}^{(j)})$ is the full
$\sigma$-algebra of the finite copy, the precise definition is
\begin{equation}\label{eq:TH-filtration-recursion}
 \begin{aligned}
 \mathcal G^{(N)}_{i,\alpha}
 := {}&\operatorname{pr}_0^{-1}\mathcal A_i
   \vee\bigvee_{j\prec i}
            \operatorname{pr}_j^{-1}\mathscr P(W_{N-1}^{(j)})
 \vee\operatorname{pr}_i^{-1}\mathcal G^{(N-1)}_\alpha.
 \end{aligned}
\end{equation}
The symbol $\vee$ means the $\sigma$-algebra generated by the indicated
$\sigma$-algebras. At level zero the only $\sigma$-algebra is the trivial
one on $W_0$.

This is an increasing list of $\sigma$-algebras.
Within the block belonging to $i$, the earlier copies and the
root information remain fixed, while the information in the
current copy increases by induction.
When we pass from one child block to the next, all information
in the preceding copy is included among the fully revealed
earlier copies. At that transition the root information also
increases or stays unchanged, since the list
\eqref{eq:TH-seed-filtration} is increasing.
The first $\sigma$-algebra is trivial by induction and the definition
of the first $\mathcal A_i$.
At the last index the root sign is known, every preceding copy is
fully known, and the last copy has reached its full $\sigma$-algebra.
Thus the last $\sigma$-algebra is the full one on $W_N$.
This describes explicitly the filtration arrangement used in
\cite[Theorem~3.5]{Qiu}.

We check that the list is allowed in our finite-cube definition
of the Stein constants.
Let $M_N$ be the number of independent signs in $W_N$.
The product construction \eqref{eq:TH-branch-spaces} gives
$M_0=0$ and $M_N=1+4M_{N-1}$.
Iterating gives
$M_N=1+4+\cdots+4^{N-1}=(4^N-1)/3$.
Thus $W_N$ is a uniform cube of dimension $M_N$.
By \eqref{eq:TH-filtration-recursion}, every $\sigma$-algebra in the
list is generated by a specified subset of these signs.
The subsets increase with the index, and the last contains all
$M_N$ signs.
For each sign, record the first index at which it belongs to
the subset, and order the signs by these indices, with arbitrary
ordering among ties.
At any index $\alpha$, its subset is now precisely the first
$k_\alpha$ signs, where $k_\alpha$ is the cardinality of that
subset. After this permutation of cube coordinates,
\begin{equation}\label{eq:TH-canonical-filtration}
 \mathcal G^{(N)}_\alpha
 =\mathscr F_{k_\alpha}^{(M_N)},
 \qquad 0\leq k_\alpha\leq M_N.
\end{equation}
Permuting the coordinates is a probability-preserving bijection.
Composing functions with this bijection preserves their $L^2$-norms
and carries each conditional expectation to the corresponding
canonical one, since both are averages over the unrevealed signs.
The list has $|I^N|=4^N$ indices. Repeated values of $k_\alpha$
and jumps revealing several signs at once are permitted in
\eqref{eq:PX-finite-Stein}.
Hence the entire list gives admissible tests with $4^N$ summands.

We define the linear operator $T_N$ on functions $f:W_N\to H_N$ by
\begin{equation}\label{eq:TH-coordinate-expectations}
 (T_Nf)_{i,\alpha}
 :=\E[f_{i,\alpha}\mid\mathcal G^{(N)}_{i,\alpha}].
\end{equation}
Each scalar coordinate uses its assigned conditional expectation.
Inductively, $(h_{N-1})_\alpha$ is measurable with respect to
$\mathcal G^{(N-1)}_\alpha$.
Also $\theta_i$ is $\mathcal A_i$-measurable by
\eqref{eq:TH-seed} and \eqref{eq:TH-seed-filtration}.
Both factors in \eqref{eq:TH-h-recursion} are
therefore measurable with respect to
$\mathcal G^{(N)}_{i,\alpha}$, so the conditional expectation
fixes their product. By \eqref{eq:TH-coordinate-expectations},
coordinate by coordinate this gives
\begin{equation}\label{eq:TH-fixed-point}
 T_Nh_N=h_N.
\end{equation}

\subsection{Accumulation of the centralizer residual}

We first record how the centralizer acts on the product form used
in \eqref{eq:TH-h-recursion}.
Let $\theta\in H_1$ and $g^{(i)}\in H_{N-1}$ satisfy
$\|\theta\|_{H_1}=\|g^{(i)}\|_{H_{N-1}}=1$, and put
$h_{i,\alpha}:=\theta_i g^{(i)}_\alpha$.
With $\Omega_0=0$ on $H_0=\mathbb C$, we claim that
\begin{equation}\label{eq:TH-Omega-recursion}
 (\Omega_N(h))_{i,\alpha}
 =(\Omega_1(\theta))_i g^{(i)}_\alpha
       +\theta_i(\Omega_{N-1}(g^{(i)}))_\alpha.
\end{equation}
To compute the first averages, write $i=(u,v)$.
Averaging over all the child coordinates gives
$|\theta_i|^2\|g^{(i)}\|_{H_{N-1}}^2=|\theta_i|^2$.
Consequently
\begin{align*}
 a_1(h)&=1=a_1(\theta),\\
 a_2(h)(u)&=\E_v|\theta_{u,v}|^2=a_2(\theta)(u),\\
 a_3(h)(i)&=|\theta_i|^2=a_3(\theta)(i).
\end{align*}
At a nonzero coordinate of $h$, the first two logarithms in
\eqref{eq:TH-Omega} are thus exactly those defining
$(\Omega_1(\theta))_i$, multiplied by $g^{(i)}_\alpha$.
For $N\geq2$, the remaining averages are
\[
 a_{q+2}(h)(i,\cdot)
 =|\theta_i|^2 a_q(g^{(i)})(\cdot)
 \qquad(1\leq q\leq2N-1).
\]
In their consecutive ratios the common factor $|\theta_i|^2$
cancels. The signs agree with the child's signs because
$(-1)^{q+3}=(-1)^{q+1}$.
The remaining $2N-2$ logarithms therefore give
$\theta_i(\Omega_{N-1}(g^{(i)}))_\alpha$.
If $\theta_i=0$, both terms on the right of
\eqref{eq:TH-Omega-recursion} are zero by the zero-coordinate
convention. If $g^{(i)}_\alpha=0$, the same is true.
This proves the formula at every coordinate.
At $N=1$ there are no child variables, $g^{(i)}\in H_0$
is a scalar of modulus one, and there are no further logarithms.
The same formula follows with $\Omega_0=0$.

For the functions already constructed, define
\begin{equation}\label{eq:TH-C-definition}
 C_N:=T_N\Omega_N(h_N)-\Omega_N(h_N),\qquad C_0:=0.
\end{equation}
We now compute its recursion.
The root sign, the current child space, and all other child
spaces are independent factors of the product probability
in \eqref{eq:TH-branch-spaces}.
For scalar functions $a$ of the root and $b$ of the current child,
this gives
\begin{equation}\label{eq:TH-conditional-factorisation}
 \begin{aligned}
 \E[a(w_0)b(w_i)\mid\mathcal G^{(N)}_{i,\alpha}](w)
 = \E[a\mid\mathcal A_i](w_0)\,
 \E[b\mid\mathcal G^{(N-1)}_\alpha](w_i).
 \end{aligned}
\end{equation}
To verify it explicitly, let $B_0$ be the atom of $\mathcal A_i$
containing $w_0$, and let $B_i$ be the atom of the child's
$\mathcal G^{(N-1)}_\alpha$ containing $w_i$.
An atom of \eqref{eq:TH-filtration-recursion} also fixes all
earlier children and leaves the later children unrestricted.
Averaging a function depending only on $w_0,w_i$ over this atom
therefore reduces to
\[
 \frac1{|B_0|\,|B_i|}
 \sum_{s\in B_0}\sum_{t\in B_i}a(s)b(t)
 =\Bigl(\frac1{|B_0|}\sum_{s\in B_0}a(s)\Bigr)
  \Bigl(\frac1{|B_i|}\sum_{t\in B_i}b(t)\Bigr).
\]
These are the two conditional expectations in
\eqref{eq:TH-conditional-factorisation}.

We apply \eqref{eq:TH-Omega-recursion} to the seed $\theta(w_0)$
and the child vectors $h_{N-1}(w_i)$, which are normalised by
\eqref{eq:TH-h-normalisation}.
For the first term of that formula, the child's coordinate is
fixed by its designated conditional expectation, by
\eqref{eq:TH-fixed-point} at level $N-1$.
For the second term, $\theta_i$ is $\mathcal A_i$-measurable.
Thus \eqref{eq:TH-conditional-factorisation} gives
\begin{align*}
 (T_N\Omega_N(h_N))_{i,\alpha}
 &=\E[(\Omega_1(\theta))_i\mid\mathcal A_i](w_0)
                (h_{N-1}(w_i))_\alpha\\
 &\qquad\qquad +\theta_i(w_0)
        (T_{N-1}\Omega_{N-1}(h_{N-1})(w_i))_\alpha.
\end{align*}
Subtracting \eqref{eq:TH-Omega-recursion} and using
\eqref{eq:TH-seed-residual-definition} for $c$ and
\eqref{eq:TH-C-definition} for $C_{N-1}$ proves
\begin{equation}\label{eq:TH-C-recursion}
 (C_N(w))_{i,\alpha}
 =c_i(w_0)(h_{N-1}(w_i))_\alpha
       +\theta_i(w_0)(C_{N-1}(w_i))_\alpha.
\end{equation}

We prove by induction the identities
\begin{equation}\label{eq:TH-branch-identities}
 \begin{aligned}
 \|h_N(w)\|_{H_N}^2&=1,\\
 \langle h_N(w),C_N(w)\rangle_{H_N}&=0,\\
 \|C_N(w)\|_{H_N}^2&=N\kappa^2,
 \qquad w\in W_N.
 \end{aligned}
\end{equation}
They hold at level zero because $h_0=1$ and $C_0=0$.
Suppose that they hold at level $N-1$ and fix $w\in W_N$.
For this calculation abbreviate
$h_i:=h_{N-1}(w_i)$ and $C_i:=C_{N-1}(w_i)$,
and evaluate $\theta_i,c_i$ at $w_0$.
The product Hilbert norm uses a factor $1/4$ for the root
coordinate. By \eqref{eq:TH-h-recursion} and
\eqref{eq:TH-C-recursion},
\[
 \langle h_N,C_N\rangle_{H_N}
 =\frac14\sum_{i\in I}
     \bigl(\theta_i c_i\|h_i\|_{H_{N-1}}^2
          +\theta_i^2\langle h_i,C_i\rangle_{H_{N-1}}\bigr).
\]
All coordinates are real, by construction.
The induction hypothesis gives $\|h_i\|^2=1$ and
$\langle h_i,C_i\rangle=0$ for each child.
The last display reduces to
$\frac14\sum_i\theta_i c_i=\langle\theta,c\rangle_{H_1}=0$
by \eqref{eq:TH-seed-identities}.

For the squared norm, the same recursion gives
\begin{align*}
 \|C_N\|_{H_N}^2
 =\frac14\sum_{i\in I}\bigl(
 &|c_i|^2\|h_i\|_{H_{N-1}}^2
 +|\theta_i|^2\|C_i\|_{H_{N-1}}^2
 +2c_i\theta_i\langle h_i,C_i\rangle_{H_{N-1}}\bigr).
\end{align*}
Each mixed term vanishes separately by the child's orthogonality.
The other two terms sum to
\[
 \frac14\sum_i|c_i|^2
  +(N-1)\kappa^2\,\frac14\sum_i|\theta_i|^2
 =\kappa^2+(N-1)\kappa^2=N\kappa^2,
\]
using \eqref{eq:TH-seed-identities} again.
The normalisation of $h_N$ was proved in
\eqref{eq:TH-h-normalisation}, completing the induction.
In particular, $\|C_N(w)\|_{H_N}=\kappa\sqrt N$ for every
$w\in W_N$. The lower bound in the final Stein test will therefore
hold before any expectation over $W_N$ is taken.

\subsection{The Stein lower bound}

\begin{proof}[Proof of Theorem~\ref{thm:TH-main} and
Corollary~\ref{cor:TH-three-space}]
Let $N\geq1$. For $\alpha\in I^N$, let $e_\alpha$ denote the
coordinate indicator on $I^N$, and let $r_\alpha$ be independent
Rademacher variables on a probability space independent of $W_N$.
We define the $Z_N$-valued functions
\begin{equation}\label{eq:TH-test-inputs}
 f_\alpha(w):=
 \bigl((\Omega_N(h_N(w)))_\alpha e_\alpha,\,
                    (h_N(w))_\alpha e_\alpha\bigr).
\end{equation}
Since $W_N$ is finite, these are measurable $L^2$ inputs.
Each input has both components supported at the single coordinate
$\alpha$, although its scalar coefficients can depend on all of $w$.

For each fixed choice of the external signs, put
\begin{equation}\label{eq:TH-diagonal-signs}
 D_rx:=(r_\alpha x_\alpha)_{\alpha\in I^N}.
\end{equation}
Multiplication by these signs preserves $|x_\alpha|^2$ at every
coordinate and hence every partial average in
\eqref{eq:TH-square-averages}.
The formula \eqref{eq:TH-Omega}, including its zero-coordinate
convention, therefore gives
\begin{equation}\label{eq:TH-sign-equivariance}
 \Omega_N(D_rx)=D_r\Omega_N(x).
\end{equation}
Also $D_r$ from \eqref{eq:TH-diagonal-signs} is an isometry on
$H_N$ by \eqref{eq:TH-Hilbert-norm}.
It follows from \eqref{eq:TH-qnorm} that
\[
 q_N(D_ru,D_rx)
 =\|D_r(u-\Omega_N(x))\|_{H_N}+\|D_rx\|_{H_N}
 =q_N(u,x).
\]
Summing the inputs in \eqref{eq:TH-test-inputs} gives
\[
 \sum_{\alpha\in I^N}r_\alpha f_\alpha(w)
 =\bigl(D_r\Omega_N(h_N(w)),D_rh_N(w)\bigr).
\]
By \eqref{eq:TH-sign-equivariance}, the first component is
$\Omega_N(D_rh_N(w))$. The second identity in
\eqref{eq:TH-block-identities} therefore gives the norm of the
pair as
\[\|D_rh_N(w)\|_{H_N}=\|h_N(w)\|_{H_N}=1,\]
by \eqref{eq:TH-h-normalisation}. Hence
\begin{equation}\label{eq:TH-test-input-bound}
 \Bigl\|\sum_\alpha r_\alpha f_\alpha(w)\Bigr\|_{Z_N}=1,
 \qquad w\in W_N.
\end{equation}

We apply conditional expectation with respect to each
$\mathcal G^{(N)}_\alpha$ to the corresponding input.
The definition \eqref{eq:TH-coordinate-expectations} of $T_N$
and the fixed-point identity
\eqref{eq:TH-fixed-point} give
\[
 \sum_\alpha r_\alpha
      \E[f_\alpha\mid\mathcal G^{(N)}_\alpha](w)
 =\bigl(D_rT_N\Omega_N(h_N)(w),D_rh_N(w)\bigr).
\]
Its first residual equals
\begin{align*}
 &D_rT_N\Omega_N(h_N)(w)-\Omega_N(D_rh_N(w))\\
 &\qquad=D_r\bigl(T_N\Omega_N(h_N)(w)-\Omega_N(h_N(w))\bigr)
 =D_rC_N(w),
\end{align*}
by \eqref{eq:TH-sign-equivariance} and
\eqref{eq:TH-C-definition}.
Thus \eqref{eq:TH-branch-identities} gives
\[
 q_N\Bigl(\sum_\alpha r_\alpha
      \E[f_\alpha\mid\mathcal G^{(N)}_\alpha](w)\Bigr)
 =\|C_N(w)\|_{H_N}+\|h_N(w)\|_{H_N}
 =1+\kappa\sqrt N.
\]
The lower bound in \eqref{eq:TH-uniform-equivalence} now gives
\begin{equation}\label{eq:TH-test-output-bound}
 \Bigl\|\sum_\alpha r_\alpha
      \E[f_\alpha\mid\mathcal G^{(N)}_\alpha](w)\Bigr\|_{Z_N}
 \geq b(1+\kappa\sqrt N).
\end{equation}
Both \eqref{eq:TH-test-input-bound} and
\eqref{eq:TH-test-output-bound} hold for every $w$ and every
choice of the external signs.
Squaring and integrating over the two probability spaces gives
an input $L^2$-norm equal to one and an output $L^2$-norm at least
$b(1+\kappa\sqrt N)$.
By \eqref{eq:TH-canonical-filtration}, these conditional
expectations form a test of \eqref{eq:PX-finite-Stein} on a
uniform cube, with $4^N$ summands.
The definition \eqref{eq:TH-finite-Stein-constants} of the finite
constant therefore gives
\begin{equation}\label{eq:TH-block-Stein}
 \mathfrak s_{4^N}(Z_N)\geq b(1+\kappa\sqrt N).
\end{equation}

By the direct sum definition \eqref{eq:TH-sums}, the coordinate
inclusion $Z_N\to Z$ is isometric.
Composing each $f_\alpha$ with that inclusion leaves both norms
in the test unchanged, and the inclusion commutes with the finite
conditional averages. Hence \eqref{eq:TH-block-Stein} also gives
the same lower bound for $\mathfrak s_{4^N}(Z)$.
For an integer $n\geq4$, put $N:=\lfloor\log_4 n\rfloor$.
Then $4^N\leq n<4^{N+1}$ and $N\geq1$.
Appending $n-4^N$ zero inputs produces a test with $n$ summands,
so
\[
 \mathfrak s_n(Z)\geq b(1+\kappa\sqrt N).
\]
Since $n$ is an integer, $n+1\leq4^{N+1}$, and consequently
$\log(n+1)\leq(N+1)\log4\leq2N\log4$.
It follows that
\[
 \mathfrak s_n(Z)
 \geq\frac{b\kappa}{\sqrt{2\log4}}\sqrt{\log(n+1)},
 \qquad n\geq 4.
\]
For $n=2,3$ we use $\mathfrak s_n(Z)\geq1$ from
\eqref{eq:TH-finite-Stein-basic}.
Thus \eqref{eq:TH-main-bounds} holds, for
example, with
$c=\min\{b\kappa/\sqrt{2\log4},\,1/\sqrt{\log4}\}$.
Here $b$ and $\kappa$ are the numerical constants in
\eqref{eq:TH-uniform-equivalence} and \eqref{eq:TH-seed-kappa}. In particular,
\[\sup_n\mathfrak s_n(Z)=\infty,\]
so \eqref{eq:TH-finite-Stein-criterion} shows that $Z$ fails the Stein property.
The exact sequence and superreflexivity of $Z$ were established in
Subsection~\ref{subsec:TH-exact-sequence}, so this completes the
proof of Theorem~\ref{thm:TH-main}.

Finally, every Hilbert space has the Stein property.
In fact, for Hilbert-valued inputs on a finite cube, expansion
of the squared norm and $\E_r(r_ir_j)=\mathbf1_{\{i=j\}}$ give
\[
 \E_r\Bigl\|\sum_j r_jf_j\Bigr\|_{L^2(H)}^2
 =\sum_j\|f_j\|_{L^2(H)}^2.
\]
The same identity applies to the conditional expectations of
the inputs. Since each conditional expectation is an $L^2(H)$
contraction, summing its squared contraction bounds proves
\eqref{eq:PX-finite-Stein} with constant one.
The property is also invariant under Banach space isomorphisms.
Indeed, a bounded linear isomorphism commutes with the finite
conditional averages. Applying the Stein inequality after the
isomorphism, and then applying its inverse to the output,
multiplies the constant by at most the product of their norms.
By \eqref{eq:TH-kernel-embedding}, $j(H)$ is isomorphic to $H$,
and \eqref{eq:TH-quotient-isometry} shows that $Z/j(H)$ is isometric to $H$.
Both therefore have the Stein property, whereas $Z$ fails it.
This proves Corollary~\ref{cor:TH-three-space}.
\end{proof}

\bigskip

\

\bigskip

\subsection*{Acknowlegdment} I thank Mark Veraar for the suggestion to look into this problem.

\subsection*{AI disclosure statement}

This manuscript gradually emerged as a result of long conversations with OpenAI's ChatGPT, covering exploratory mathematical discussion, checking arguments, and polishing the text. The author has checked all arguments and takes full responsibility for the mathematical
content and the final version of the paper. The Lean formalisation of Theorem \ref{thm:TH-main} was developed with the assistance of OpenAI's Codex.

\end{document}